# GENEALOGY OF CATALYTIC BRANCHING MODELS[1]

By Andreas Greven, Lea Popovic and Anita Winter

*University Erlangen–Nurnberg, Concordia University and
University Erlangen–Nurnberg*

We consider catalytic branching populations. They consist of a catalyst population evolving according to a critical binary branching process in continuous time with a constant branching rate and a reactant population with a branching rate proportional to the number of catalyst individuals alive. The reactant forms a process in random medium.

We describe asymptotically the genealogy of catalytic branching populations coded as the induced forest of $\mathbb{R}$-trees using the many individuals—rapid branching continuum limit. The limiting continuum genealogical forests are then studied in detail from both the quenched and annealed points of view. The result is obtained by constructing a contour process and analyzing the appropriately rescaled version and its limit. The genealogy of the limiting forest is described by a point process. We compare geometric properties and statistics of the reactant limit forest with those of the "classical" forest.

**1. Introduction.** It is well known that Feller's branching diffusion and, more generally, continuous state branching processes appear as weak limits of a broad range of suitably rescaled Galton–Watson processes (see, e.g., [18]). Provided one introduces a suitable topology, the rescaled family forests of these Galton–Watson processes converge to limit forests which can be represented by paths of a reflected Brownian motion, and more generally a reflected Lévy process without negative jumps (see [1, 2, 3, 11, 19] and references therein). These invariance principles are the main tool for analyzing the asymptotic behavior of genealogies of Galton–Watson processes. For example, one can conclude that the height of a Galton–Watson tree with total population size $n$ behaves in distribution, as $n \to \infty$, as $2\sqrt{n}$ times the height of a standard Brownian excursion.

Received June 2006; revised October 2008.

[1]Supported by DFG-Forschergruppe 498 through Grant GR 876/13-1.

One of the next important steps is to investigate genealogies of multi-type branching models with possible interactions between populations. For neutral two-type branching models there are three universality classes of behavior: independent branching, (one-sided) catalytic branching and mutually catalytic branching. Loss of independence in the two latter models generates many new features and causes great technical difficulties for its analysis.

In the present paper we focus on describing the genealogy of a catalytic branching diffusion. This is the many individual fast branching limit of an interacting branching particle model involving two populations, in which one population, the "catalyst," evolves autonomously according to a Galton–Watson process, while the other population, the "reactant," evolves according to branching dynamics that are dependent on the number of catalyst particles.

We show that the sequence of suitably rescaled family forests for the reactant population converges in Gromov–Hausdorff topology to a limit forest and we characterize the distribution of this limit. We first describe the limit of the reactant family forest cut off at a height where the catalyst falls below a certain threshold. We show that its contour process is given by a path of a reflected diffusion process. From that path we derive a collection of point processes each describing the mutual genealogical distances between all individuals in the population alive at a certain time.

We start the analysis from a "quenched" point of view, that is, we first fix a realization of the catalyst population as an autonomous branching process, and consider the reactant population as a branching process in a given medium. We then change to an "annealed" point of view and show that the joint law of the rescaled catalyst and the reactant population converges to a limit.

Finally, we discuss differences between the reactant limit forest and the classical continuum random forest which is known to be associated with a reflected Brownian motion. The point processes exhibit the difference in genealogical structure between the reactant population family forests and the classical Brownian continuum random tree. The differences are due to two different effects: the vanishing branching rate once the catalyst becomes extinct and the temporal inhomogeneity of the branching rates. The differences are roughly as follows. If we pick from the population of the reactant at time $t$ two individuals at random, their chance of being from different trees is smaller than in the classical case (with a suitable choice of the branching rate). Furthermore the total tree length becomes infinite in a stronger sense, namely instead of having Hausdorff dimension 2 with finite Hausdorff measure function, we now have a Hausdorff measure function that grows faster than a square. As a consequence, with positive probability the reactant limit forest cannot be associated with the path of a diffusion.

*Outline.* The rest of the paper is structured as follows. In Section 2 we provide the underlying spaces and the topologies we will need for representing the genealogical forests of our branching model. In Section 3 we define our catalytic particle model, its genealogical forests, the associated contour processes, the genealogical point processes and the suitably rescaled versions of all of these objects. In Section 4 we state our results for the reactant limit forest and its differences from the classical Brownian forest. Subsequent sections contain proofs of these results. In Section 5 we show that the family of suitably rescaled catalyst forests has a limit forest. In Section 6 we prove the statement on the limit of the reactant contour processes. In Section 7 we use the limit contour process to uniquely characterize the distribution of the limit forest. In Section 8 we prove results on the limit reactant genealogical point process given the catalyst medium. In Section 9 we show that the joint law of catalyst and reactant tree converges. Finally Section 10 contains calculations establishing the differences between the reactant and the classical forest.

## 2. Preliminaries on trees and their representations.

In this section we introduce the relevant spaces and topologies for both the particle model and its diffusion limit.

The standard way of representing genealogical relationships between individuals in a branching population is with a family forest. The family forest consists of as many trees as the initial number of individuals. Each individual in a branching process has an edge associated to it whose length is equal to its lifetime. When an individual gives birth the edge branches into two new edges, while its death turns its edge end into a leaf. In particular, time in the branching process corresponds to height in the family forest which is measured in terms of the distance to the roots.

A family forest can be embedded in a plane. A planar embedding implicitly requires a linear (total) ordering to be placed on individuals, providing each with a label. A lexicographic labeling, also called the Ulam–Harris labeling, gives each individual a label from the set $\mathbb{N}^n$, where $n \in \{1, 2, \ldots\}$ is the number of ancestors the individual has had from the start of the process. The initial individuals are given labels between 1 and the initial population size, distributed in a random order to them. Each parent subsequently gives a label to its children that consists of its own label followed by a number between 1 and the total number of its children, distributing these randomly among its children. See Figure 2 for more detail, and [21] or [11] for a formal definition.

A planar representation and a lexicographic order are the most convenient way to visualize a family forest with a linear order. However, the appropriate space for family trees in general is the set of all rooted $\mathbb{R}$-trees with a linear order on them.

*Rooted linearly ordered $\mathbb{R}$-trees.* A metric space $(T, d)$ is called an $\mathbb{R}$-*tree* if it is path-connected and satisfies the so-called "four point condition," that is, $d(x_1, x_2) + d(x_3, x_4) \leq \max\{d(x_1, x_3) + d(x_2, x_4), d(x_1, x_4) + d(x_2, x_3)\}$, for all $x_1, \ldots, x_4 \in T$. We also refer the reader to [6, 8, 9, 10, 27] for a background on $\mathbb{R}$-trees. A *rooted $\mathbb{R}$-tree*, $(T, d, \rho)$, is an $\mathbb{R}$-tree $(T, d)$ with a distinguished point $\rho \in T$ that we call the *root*. A forest is just a collection of rooted $\mathbb{R}$-trees, and we can always associate a forest of rooted $\mathbb{R}$-trees with an $\mathbb{R}$-tree with the root being of degree greater than 1 by gluing together all the roots. In our paper we will think of a forest as a rooted $\mathbb{R}$-tree whose root will have degree 1 or larger.

Let $\mathbb{T}^{\mathrm{root}}$ denote the collection of all root-invariant isometry classes of rooted compact $\mathbb{R}$-trees, where we define a root-invariant isometry to be an isometry $i : (T_1, d_{T_1}, \rho_{T_1}) \to (T_2, d_{T_2}, \rho_{T_2})$ which in addition satisfies $i(\rho_{T_1}) = \rho_{T_2}$. Let $T = \{\iota\}$ denote the points of a metric space $(T, d, \rho) \in \mathbb{T}^{\mathrm{root}}$. For all $t \geq 0$, let $\partial Q_t$ denote the set of individuals alive at time $t$ in the branching population. Define the genealogical distance metric for $(T, d, \rho)$ by setting for any $t_1, t_2 \geq 0$ and any $\iota_1 \in \partial Q_{t_1}, \iota_2 \in \partial Q_{t_2}$

$$d(\iota_1, \iota_2) := t_1 + t_2 - 2\tau(\iota_1, \iota_2),$$

where $\tau(\iota_1, \iota_2)$ denotes the death of, or splitting time from the "most recent common ancestor" of the individuals $\iota_1$ and $\iota_2$ alive at times $t_1$ and $t_2$, respectively. Moreover, for any two $\rho, \rho' \in \partial Q_0$ by our convention we have $d(\rho, \rho') = 0$.

For an $\mathbb{R}$-tree $(T, d)$ the root $\rho \in T$ defines a partial order $\leq$ on $T$. For all $y \in T$, let $[\rho, y]$ denote the unique closed *arc* between them, that is, the image of the unique geodesic between $\rho$ and $y$. We define $x \leq y$ if $x \in [\rho, y]$. A linear order (total order) relation $\leq^{\mathrm{lin}}$ on a rooted ordered $\mathbb{R}$-tree is a total order which is compatible with the partial order [condition (i)] and fulfills the following [condition (ii); see Figure 1]:

  (i)  For all $x, y \in T$ if $x \leq y$, then $x \leq^{\mathrm{lin}} y$.
  (ii) For all $x, y \in T$ with $x \leq^{\mathrm{lin}} y$, if $x', y' \in T$ are such that $x \wedge y \leq x' \wedge x$ and $x \wedge y \leq y' \wedge y$ then $x' \leq^{\mathrm{lin}} y'$.

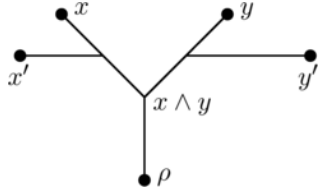

FIG. 1. *Illustrates condition* (ii): *if* $x \leq^{\mathrm{lin}} y$, $x'$ *is in the same "subtree above"* $x \wedge y$ *as* $x$, *and* $y'$ *is in the same "subtree above"* $x \wedge y$ *as* $y$, *then* $x' \leq^{\mathrm{lin}} y'$.

Then $(T, d, \rho, \leq^{\mathrm{lin}})$ is a *rooted $\mathbb{R}$-tree with a linear order*. Let $\mathbb{T}^{\mathrm{root,lin}}$ denote the set of all root-and-linear order invariant isometry classes of compact $\mathbb{R}$-trees, where we define a root-and-linear root order invariant isometry $i: (T, d_T, \rho_T, \leq_T^{\mathrm{lin}}) \to (T', d_{T'}, \rho_{T'}, \leq_{T'}^{\mathrm{lin}})$ to be a root invariant isometry $i: (T, d_T, \rho_T) \to (T', d_{T'}, \rho_{T'})$ which in addition satisfies $i(u) \leq_{T'}^{\mathrm{lin}} i(u')$ iff $u \leq_T^{\mathrm{lin}} u'$.

REMARK 2.1. For any element of $\mathbb{T}^{\mathrm{root,lin}}$ there is a planar representative of this isometry class. It is visually convenient to use it as our canonical representative when describing a family forest.

*Contour processes.* There is a long tradition of coding $\mathbb{R}$-trees via excursions (see [3, 11, 19]) which we briefly recall. Write $C^0_{\mathbb{R}^+}[0, \infty)$ for the space of continuous functions $e: [0, \infty) \to \mathbb{R}^+$ with $e(0) = 0$. We associate each $e \in C^0_{\mathbb{R}^+}[0, \infty)$ with an $\mathbb{R}$-tree as follows. Define an equivalence relation $\sim_e$ on $\mathbb{R}^+$ by letting $u_1 \sim_e u_2$ iff $e(u_1) = \inf_{u \in [u_1 \wedge u_2, u_1 \vee u_2]} e(u) = e(u_2)$. Consider the following pseudo-metric on $[0, 1] : d_{T_e}(u_1, u_2) := e(u_1) + e(u_2) - 2 \inf_{u \in [u_1 \wedge u_2, u_1 \vee u_2]} e(u)$, which becomes a true metric on the quotient space $T_e := \mathbb{R}_+|_{\sim_e}$. Then by Lemma 3.1 in [14], for each $e \in C^0_{\mathbb{R}^+}[0, \infty)$ with $\sup\{t \geq 0 : e(t) > 0\} < \infty$ the metric space $(T_e, d_{T_e})$ is a compact ordered $\mathbb{R}$-tree.

For $e \in C^0_{\mathbb{R}^+}[0, \infty)$, let $[0] := \{u \in [0, \infty) : e(u) = 0\}$ and define the map

$$(2.1) \qquad \mathcal{T}: C^0_{\mathbb{R}^+}[0, \infty) \to \mathbb{T}^{\mathrm{root,lin}},$$

which maps an excursion to a rooted linearly ordered $\mathbb{R}$-tree by $\mathcal{T}(e) := (T_e, d_{T_e}, [0], \leq_e^{\mathrm{lin}})$, with the linear order $\leq_e^{\mathrm{lin}}$ defined by $u_1 \leq_e^{\mathrm{lin}} u_2$ iff $\inf[u_1] \leq \inf[u_2]$, here $[u]$ is the set of elements equivalent to $u$. The map $\mathcal{T}$ is continuous (Lemma 3.1, [14]).

A rooted linearly ordered tree is *finite* if it has finitely many branch points and leaves. Let $\mathbb{T}^{\mathrm{root,lin}}_{\mathrm{fin}}$ denote the subspace in $\mathbb{T}^{\mathrm{root,lin}}$ of finite rooted linearly ordered trees, and fix a speed $\sigma > 0$. The linear order allows one to define a depth-first search on the set of its leaves and branch points. Depth-first search is a function that starting from the root visits every point of this tree while respecting its linear order by traversing the tree in such a way that its first visit to $x$ comes before its first visit to $y$ if $x \leq^{\mathrm{lin}} y$. If one associates a speed $\sigma$ with this traversal and assigns a function which records the height of the tree at every time point of the traversal, then as a result one obtains a function $e \in C^0_{\mathbb{R}^+}[0, \infty)$.

This defines a *contour process* map $\mathcal{C}(\cdot; \sigma): \mathbb{T}^{\mathrm{root,lin}}_{\mathrm{fin}} \to C^0_{\mathbb{R}^+}[0, \infty)$, which associates with each $T \in \mathbb{T}^{\mathrm{root,lin}}_{\mathrm{fin}}$ an excursion $\mathcal{C}(T; \sigma) \in C^0_{\mathbb{R}^+}[0, \infty)$. It is easiest to see the correspondence of a finite rooted linearly ordered tree and its contour process from its planar representative (see Figure 2). Typically, one sets the speed $\sigma = 2$, and for all $\sigma$ we have $\mathcal{T} \circ \mathcal{C} = id, \mathcal{C} \circ \mathcal{T} = id$. A formal description of the contour process can be found in Section 6.1 in [21].

*Genealogical point-processes.* We next consider finite rooted linearly ordered trees spanned by a population alive at a fixed time. Such trees are *balanced*, that is, have all the leaves at the same height. Given a family tree of a branching process and a fixed time $t > 0$ this is the subtree of the family tree containing only information about the genealogy of the individuals alive at time $t$. That is, it is defined by the collection of mutual distances between the individuals alive at time $t$. This collection of distances forms an ultra-metric space, that is, for all $x_1, x_2, x_3$ with $d(\rho, x_1) = d(\rho, x_2) = d(\rho, x_3)$ we have $d(x_1, x_3) \leq d(x_1, x_2) \vee d(x_2, x_3)$.

For a fixed $t > 0$ we can record these distances in a point process. The linear order gives a linear order on the individuals alive at time $t$. We can represent these individuals as points along a line at height $t$ with equal spacing between them. We can record the mutual distance between each neighboring pair of individuals as a point equally distant from the points in this pair at a height smaller than $t$. More formally, given $t > 0$ and $(T, \rho) \in \mathbb{T}_{\text{fin}}^{\text{root,lin}}$, let

$$(2.2) \qquad Q_t(T, \rho) := \{x \in T : d(x, \rho) \leq t\}$$

and

$$\partial Q_t(T, \rho) := \{x \in T : d(x, \rho) = t\}.$$

Define a map $\mathcal{P}^t(\cdot; \varsigma)$ which sends an element in $\mathbb{T}_{\text{fin}}^{\text{root,lin}}$ to a point process on $[0, \infty) \times [0, \infty)$ as follows. Given $(T, \rho) \in \mathbb{T}_{\text{fin}}^{\text{root,lin}}$, $\partial Q_t(T, \rho)$ is an ordered finite set, so we can write $\partial Q_t(T, \rho) := \{x_1, \ldots, x_{\#\partial Q_t(T, \rho)}\}$ assuming that $x_1 \leq^{\text{lin}} x_2 \leq^{\text{lin}} \cdots \leq^{\text{lin}} x_{\#\partial Q_t(T, \rho)}$. Let

$$\mathcal{P}^t((T, \rho); \varsigma) := \{(i\varsigma, t - \tfrac{1}{2}d(x_i, x_{i+1})); i = 1, \ldots, \#\partial Q_t(T, \rho) - 1\}.$$

The set $\mathcal{P}^t((T, \rho); \varsigma)$ is a *point-process representation* of the tree $\partial Q_t(T, \rho)$. Typically the spacing between the points is $\varsigma = 1$.

Note that mutual distances of all neighboring pairs is sufficient to reconstruct mutual distances between all pairs. For any $x_i, x_l \in \partial Q_t(T, \rho)$ such that $x_i \leq^{\mathrm{lin}} x_j$ the distance $d(x_i, x_l) = \inf_{i \leq j < k \leq l} \{d(x_j, x_k)\}$. More detailed description of the point-process representation can be found in [22].

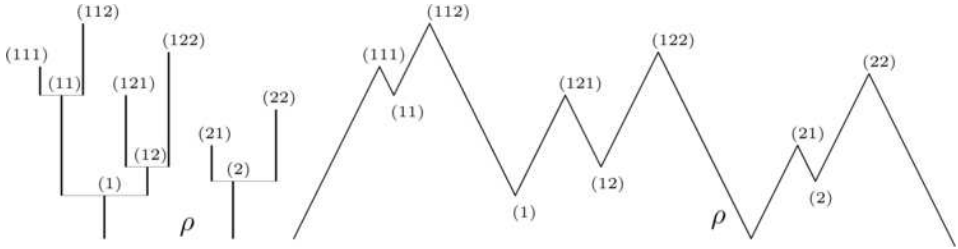

Fig. 2. *Illustrates mapping of a finite linearly ordered tree to an excursion.*

Letting the time $t > 0$ vary, this collection of subtrees of individuals alive at time $t$ defines the full family tree. That is the set $T = \bigcup_{t \geq 0} \partial Q_t(T, \rho)$ with the linear order and distance metric of $(T, d, \rho, \leq^{\mathrm{lin}})$ preserved in the subtrees spanned by $\partial Q_t(T, \rho)$. Hence, the collection of point process $\{\mathcal{P}^t((T, \rho); \varsigma)\}_{t \geq 0}$ gives a complete representation of the metric space $(T, d, \rho, \leq^{\mathrm{lin}})$.

*Gromov–Hausdorff distance.* We finally recall the appropriate measure of distance between two $\mathbb{R}$-trees. The *Hausdorff distance* between two compact subsets $K_1$ and $K_2$ of a complete metric space $(E, r)$ is defined as $d_{\mathrm{H}}(K_1, K_2) := \inf\{\varepsilon > 0 : K_1^\varepsilon \supseteq K_2, K_2^\varepsilon \supseteq K_1\}$, where for $\varepsilon > 0$ and $K \subseteq E$ closed, $K^\varepsilon := \{x \in E : r(x, K) < \varepsilon\}$. The *rooted Gromov–Hausdorff distance* between two rooted $\mathbb{R}$-trees $(T_1, \rho_{T_1})$ and $(T_2, \rho_{T_2})$ is then defined as

$$(2.3) \quad d_{\mathrm{GH^{root}}}((T_1, \rho_{T_1}), (T_2, \rho_{T_2})) := \inf d_{\mathrm{H}}^{(Z, d_Z)}(T_1', T_2') \vee d_Z(\rho_{T_1}', \rho_{T_2}'),$$

where the infimum is taken over all rooted $\mathbb{R}$-trees $(T_1', \rho_{T_1}')$ and $(T_2', \rho_{T_2}')$ that are root-invariant isomorphic to $(T_1, \rho_{T_1})$ and $(T_2, \rho_{T_2})$, respectively, and that are (as unrooted trees) subspaces of a common metric space $(Z, d_Z)$. By Theorem 2 in [13], $(\mathbb{T}^{\mathrm{root}}, d_{\mathrm{GH^{root}}})$ is a Polish space.

*Notational conventions.* Throughout the rest of the paper the symbol $\underset{n \to \infty}{\Longrightarrow}$ will always denote convergence in distribution in the underlying process spaces. These are:

- *for population size/total mass processes*: cádlág paths with Skorohod topology;
- *for forest-valued processes*: $(\mathbb{T}^{\mathrm{root}}, d_{\mathrm{GH^{root}}})$ of (isometry classes of) compact rooted $\mathbb{R}$-trees with rooted Gromov–Hausdorff topology [compare (2.3)];
- *for contour processes*: nonnegative continuous functions with the topology of uniform convergence on compact sets;
- *for point process*: locally finite measures on $[0, \infty) \times [0, \infty)$ with the weak topology (i.e., the topology given by convergence of all bounded continuous functions).

## 3. Catalytic branching model.

In this section we introduce the model and the family forests with their contour processes and point-process representation that will be used to obtain our results. We begin in Section 3.1 by defining the catalytic branching particle model and the family forests of its two populations. We define the contour process and the collection of point-processes associated with these finite family forests. We then define in Section 3.2 the suitably rescaled versions of these processes, and the diffusion limits of interest. Finally, in Section 3.3 we recall some known results for these processes which are mainly for the catalytic population.

3.1. *Definition of the model and its family forests.* The catalytic branching model consists of two different populations of distinct individuals: the *catalyst* population and the *reactant* population. The pair of populations $(\eta, \xi)$ evolves as a Markov process $(\eta_t, \xi_t)_{t \geq 0}$ according to the following rules. The catalyst population $\eta$ is a classical continuous-time critical branching process: every individual has an exponential lifetime with parameter $b_1$ after which it is replaced by 0 or 2 offspring with equal probability. The reactant population $\xi$ evolves analogously except that the exponential lifetime distribution of each individual is replaced by a lifetime distribution $F(t) = 1 - \exp(-\int_s^{s+t} b(u) \, du)$ where $s$ is the birth time of this individual and $b(u)$ equals $b_2$ times the total number of catalyst individuals at time $u$.

The total number of individuals in the catalyst, reactant population is a continuous-time Markov process with values in $(m\mathbb{N})^2$ if we associate mass $m$ with an individual. Typically one sets the mass of a single particle $m = 1$. However, to obtain a reasonable limit of the rescaled catalytic branching model we will need to scale down a single particle's mass contribution.

DEFINITION 3.1 (Total mass processes).

- The catalyst total mass process $\eta^{\mathrm{tot}} = (\eta_t^{\mathrm{tot}})_{t \geq 0}$ is a critical binary Galton–Watson process with constant branching rate $b_1$

$$\eta_t^{\mathrm{tot}} \equiv (\eta_t^{\mathrm{tot}}; m) \longmapsto \begin{cases} \eta_t^{\mathrm{tot}} + m, \\ \eta_t^{\mathrm{tot}} - m, \end{cases} \quad \text{each at rate } \frac{1}{2} b_1 \frac{\eta_t^{\mathrm{tot}}}{m},$$

- given $\eta^{\mathrm{tot}}$, the reactant total mass process $\xi^{\mathrm{tot}} = (\xi_t^{\mathrm{tot}})_{t \geq 0}$ is a critical binary Galton–Watson process with time-inhomogeneous branching rate $b_2 \eta_t^{\mathrm{tot}}$

$$\xi_t^{\mathrm{tot}} \equiv (\xi_t^{\mathrm{tot}}; m) \longmapsto \begin{cases} \xi_t^{\mathrm{tot}} + m, \\ \xi_t^{\mathrm{tot}} - m, \end{cases} \quad \text{each at rate } \frac{1}{2} b_2 \eta_t^{\mathrm{tot}} \frac{\xi_t^{\mathrm{tot}}}{m}.$$

We introduce the random family forests of the catalytic branching particle populations encoded as elements in $\mathbb{T}^{\mathrm{root}}$ as follows:

DEFINITION 3.2 (Forest-valued catalytic branching particle models). Let

$$\eta^{\mathrm{for}} := \bigg( \bigcup_{t \in \mathbb{R}^+, \iota \in \partial Q_t^\eta} \{\iota\}, d_\eta, \rho^\eta \bigg),$$

$$\xi^{\mathrm{for}} := \bigg( \bigcup_{t \in \mathbb{R}^+, \iota \in \partial Q_t^\xi} \{\iota\}, d_\xi, \rho^\xi \bigg).$$

REMARK 3.1. Recall that criticality of the catalyst process implies that $\eta^{\mathrm{tot}}$ will almost surely get absorbed at 0. Let

$$T^{0,1} := \inf\{t \geq 0 : \eta_t^{\mathrm{tot}} = 0\}.$$

Notice that after the catalyst mass process is absorbed at 0, the reactant mass process gets absorbed as well. That is, in the reactant family tree the most recent common ancestor of any two points at height greater than or equal to $T^{0,1}$ has height smaller than $T^{0,1}$, or equivalently, after $T^{0,1}$ the branches simply extend to $\infty$. In the following we impose the convention to cut off all the branches of $\xi^{\mathrm{for}}$ at height $T^{0,1}$, that is, $\xi^{\mathrm{for}}$ is the subspace of the reactant family forest consisting of all points in the reactant family forest whose distance to the root is at most $T^{0,1}$.

The forests $\eta^{\mathrm{for}}$ and $\xi^{\mathrm{for}}$ (with the convention given in Remark 3.1) have finitely many leaves and branch points, almost surely. We define their contour process representation:

DEFINITION 3.3 (Contour processes of the finite family forests). Let

$$B := (B_u)_{u \geq 0} \equiv \mathcal{C}(\eta^{\mathrm{for}}; \sigma), \qquad C := (C_u)_{u \geq 0} \equiv \mathcal{C}(\xi^{\mathrm{for}}; \sigma).$$

For the collection of genealogical point processes of the family forests in the catalytic branching model we define:

DEFINITION 3.4 (Genealogical point processes). Let

$$\Pi = \{\Pi^t\}_{t \geq 0}, \qquad \Pi^t \equiv \mathcal{P}^t(\eta^{\mathrm{for}}; \varsigma), \qquad \Xi = \{\Xi^t\}_{t \geq 0}, \qquad \Xi^t \equiv \mathcal{P}^t(\xi^{\mathrm{for}}; \varsigma).$$

3.2. *Rescaling of the model.* We now consider a family $(\tilde{\eta}^n, \tilde{\xi}^n)$ of catalytic branching particle models indexed by $n$, where the parameter $n$ means that the initial number of individuals is increased by a factor $n$ and the branching rate of the catalyst is sped up by a factor $n$.

We will define for the rescaled model the total mass process, the family forest, the contour process and the collection of point processes. We begin with the rescaled total mass processes.

DEFINITION 3.5 (Rescaled mass processes). The rescaled total mass process $(\tilde{\eta}^{\mathrm{tot},n}, \tilde{\xi}^{\mathrm{tot},n})$ is a continuous-time Markov chain with $(\tilde{\eta}_0^{\mathrm{tot},n}, \tilde{\xi}_0^{\mathrm{tot},n}) = (1,1)$ and

$$(\tilde{\eta}^{\mathrm{tot},n}, \tilde{\xi}^{\mathrm{tot},n}) \mapsto \begin{cases} \left(\tilde{\eta}^{\mathrm{tot},n} \pm \dfrac{1}{n}, \tilde{\xi}^{\mathrm{tot},n}\right), & \text{at rate } \dfrac{1}{2}n^2 b_1 \tilde{\eta}^{\mathrm{tot},n}, \\[2ex] \left(\tilde{\eta}^{\mathrm{tot},n}, \tilde{\xi}^{\mathrm{tot},n} \pm \dfrac{1}{n}\right), & \text{rate } \dfrac{1}{2}n^2 b_2 \tilde{\eta}^{\mathrm{tot},n} \tilde{\xi}^{\mathrm{tot},n}. \end{cases}$$

The rescaled population processes also have associated family forests for which we use the following notation:

DEFINITION 3.6 (Rescaled forest-valued processes). Let $\tilde{\eta}^{\mathrm{for},n}$ be the family forest of rooted $\mathbb{R}$-trees associated with $\tilde{\eta}^n$, and let $\tilde{\xi}^{\mathrm{for},n}$ be the family forest of rooted $\mathbb{R}$-trees associated with $\tilde{\xi}^n$.

Recall that we imposed the convention that the branches of $\tilde{\xi}^{\mathrm{for},n}$ are cut off above height

$$(3.1) \qquad \tilde{T}^{0,n} := \inf\{t \ge 0 : \tilde{\eta}_t^{\mathrm{tot},n} = 0\}.$$

Since the branching rate of the rescaled population process is accelerated by a factor of $n$, the edge lengths in the forest are shortened typically by a factor of $\frac{1}{n}$.

REMARK 3.2. Since we can always rescale time by a constant without any effect on the shape of a forest we shall henceforth assume that $b_1 = 1$. For many purposes, in particular if the catalyst medium is given, we also may assume that $b_2 = 1$.

With probability of order $\mathcal{O}(\frac{1}{n})$ the height of a Galton–Watson tree is of order $\mathcal{O}(n)$, and given a Galton–Watson process is still alive at a time of order $\mathcal{O}(n)$ its population size at that time is of order $\mathcal{O}(n)$. Hence, starting initially with $n$ particles and speeding up time by a factor of $n$ yields in the limit a Poisson number of trees, each having total population size of order $\mathcal{O}(n^2)$. Furthermore, since branching is sped up by a factor of $n$ all edges are of order $\mathcal{O}(\frac{1}{n})$. As a consequence, in order to find a nontrivial limit contour, in the $n$th approximation step we need to traverse the rescaled forests $\tilde{\eta}^{\mathrm{for},n}$ and $\tilde{\xi}^{\mathrm{for},n}$ at speed $\sigma = 2n$.

DEFINITION 3.7 (Rescaled contour processes). Let

$$\tilde{B}^n = (\tilde{B}_u^n)_{u \ge 0} \equiv \mathcal{C}(\tilde{\eta}^{\mathrm{for},n}; 2n), \qquad \tilde{C}^n = (\tilde{C}_u^n)_{u \ge 0} \equiv \mathcal{C}(\tilde{\xi}^{\mathrm{for},n}; 2n).$$

Furthermore, to keep the genealogical point processes representing the rescaled forests $\tilde{\eta}^{\mathrm{for},n}$ and $\tilde{\xi}^{\mathrm{for},n}$ in a compact set, we need to let the spacing be $\varsigma = \frac{1}{n}$.

DEFINITION 3.8 (Rescaled genealogical point processes). Let

$$\tilde{\Pi}^n = \{\tilde{\Pi}^{t,n}\}_{t \ge 0}, \qquad \tilde{\Pi}^{t,n} \equiv \mathcal{P}^t\left(\tilde{\eta}^{\mathrm{for},n}; \frac{1}{n}\right),$$

$$\tilde{\Xi}^n = \{\tilde{\Xi}^{t,n}\}_{t \ge 0}, \qquad \tilde{\Xi}^{t,n} \equiv \mathcal{P}^t\left(\tilde{\xi}^{\mathrm{for},n}; \frac{1}{n}\right).$$

3.3. *Limit of the mass processes and of the catalyst family forest.* Some of the rescaling limits for the introduced objects, in particular all those which are associated with the catalyst process, are well known. We next recall these results.

With standard techniques one can show that there exists a Markov process $(X, Y)$ with paths in $\mathcal{D}_{\mathbb{R}^+ \times \mathbb{R}^+}[0, \infty)$ such that

$$(3.2) \qquad (\tilde{\eta}^{\text{tot},n}, \tilde{\xi}^{\text{tot},n}) \underset{n \to \infty}{\Longrightarrow} (X, Y)$$

(see, e.g., [16]). Moreover, $(X, Y)$ is the unique strong solution to the following system of stochastic differential equations

$$(3.3) \qquad \begin{aligned} dX_t &= \sqrt{X_t}\, dW_t^X, \\ dY_t &= \sqrt{X_t Y_t}\, dW_t^Y \end{aligned}$$

with initial value $(X_0, Y_0) := (1, 1)$, and where $W^X = (W_t^X)_{t \geq 0}$ and $W^Y = (W_t^Y)_{t \geq 0}$ are two independent standard Brownian motions on the real line. We refer to $(X, Y)$ as the *catalytic Feller diffusion*.

Another useful fact about the catalytic Feller diffusion process is given by an explicit expression for the probability that $(X, Y)$ will hit the $y$-axis before the $x$-axis. Let for $\delta \geq 0$

$$(3.4) \qquad \begin{aligned} \tau^\delta &:= \inf\{t \geq 0 : X_t = \delta\}, \\ \rho^0 &:= \inf\{t \geq 0 : Y_t = 0\}. \end{aligned}$$

Then, as calculated in [20],

$$(3.5) \qquad \mathbb{P}\{\rho^0 < \tau^0\} = \left( \frac{4b_1}{b_2} \frac{Y_0}{X_0^2} + 1 \right)^{-1/2}.$$

Next we recall what is known about the genealogical forest of the catalyst. It was shown in [1, 2, 3] that a suitably rescaled family of Galton–Watson trees, conditioned to have total population size $n$, converges as $n \to \infty$ to the *Brownian continuum random tree* (*CRT*). With similar arguments one can get a result for the suitably rescaled catalyst forests; in fact, our proof here contains this as a special case. Also, all the limits of the contour process and genealogical point processes have been identified earlier, that we briefly recall.

Let $|\beta| = (|\beta|_u)_{u \geq 0}$ be a reflected Brownian motion and $\ell = (\ell_u)_{u \geq 0}$ its local time at level 0. Let $\ell^{-1}(t) = \inf\{u \geq 0 : \ell_u = t\}$ denote the inverse of the local time of $|\beta|$ at level 0. We have

$$(3.6) \qquad \tilde{B}^n \underset{n \to \infty}{\Longrightarrow} 2|\beta|_{\cdot \wedge \ell^{-1}(1)}.$$

It follows from (3.6) that there exists $X^{\text{for}} \in \mathbb{T}^{\text{root}}$ with

$$(3.7) \qquad \eta^{\text{for},n} \underset{n \to \infty}{\Longrightarrow} X^{\text{for}}.$$

Moreover, $X^{\text{for}}$ equals in distribution $\mathcal{T}(2|\beta|_{\cdot \wedge \ell^{-1}(1)})$ [recall that the map $\mathcal{T} : C^0_{\mathbb{R}^+}[0, \infty) \to \mathbb{T}^{\text{root,lin}}$ sends an excursion to a rooted linearly ordered compact $\mathbb{R}$-tree].

Note that if $X_0 = 1$ then $X_t^{\text{for}}$ consists of a Poisson$(\frac{1}{t})$ number of family trees. For the associated genealogical point process this translates into a Poisson$(\frac{1}{t}) - 1$ number of points at depth at least $t$ separating these trees in the forest. For a single tree the following is known by Theorem 5 in [22], $\{t_n\}$ is such that $t_n \underset{n \to \infty}{\longrightarrow} t$ for some $t \in \mathbb{R}^+$ then

$$\tilde{\Pi}^{t_n, n} \underset{n \to \infty}{\Longrightarrow} \pi^{\beta, t}.$$

The distribution of $\pi^{\beta, t}$ is the mixture over the law of $X_t$ of a Poisson point process whose intensity measure, for a given value of $X_t$, is

(3.8) $$\aleph^{\beta, t}(d\ell \times dh) = \mathbf{1}_{[0, X_t]}\, d\ell \otimes \mathbf{1}_{(0, t)} \frac{dh}{(t-h)^2}.$$

**4. Results.** In this section we state limiting results for the pair of family forests of rescaled catalytic branching populations. Since evolution of the reactant depends only on the evolution of the catalyst total mass process, we can analyze the conditional distribution of the reactant family forest given the values of the catalyst mass process. This is analogous to a "quenched" analysis for a process in random environment—the catalyst total mass process playing the role of the random environment.

We start Section 4.1 with the claim that the family of rescaled reactant forests given the catalyst total mass process has a limit forest. To determine its distribution, and thereby characterize a unique reactant limit forest, we first consider the environment up to the time in which the catalyst stays above an arbitrary positive threshold. We derive the limit of the rescaled and "truncated" contour processes, and we describe the reactant limit contour from this family of diffusions. Based on the latter we derive the rescaled genealogical point processes of the reactant. Section 4.2 gives the result on joint convergence of the catalyst and reactant family forests to a limiting pair of forests. Finally, in Section 4.3 we discuss the core differences between the genealogical structures of these two limiting forests.

4.1. *Reactant in a fixed catalytic background.* By (3.2), we can realize the rescaled catalyst total mass processes and the limiting Feller diffusion $X = (X_t)_{t \geq 0}$ on a common probability space such that for all $T > 0$

(4.1) $$\sup_{t \leq T} |\tilde{\eta}_t^{\text{tot}, n} - X_t| \underset{n \to \infty}{\longrightarrow} 0.$$

Hence, branching rates of the reactant processes $\tilde{\xi}^{\mathrm{tot},n}$ will be given by a sequence of functions in Skorohod space which converge uniformly on compacta to a *continuous* limit function $X$. We study the behavior of the rescaled reactant forests $\tilde{\xi}^{\mathrm{for},n}$ under assumption (4.1).

The first result states that given the catalyst medium, the family of rescaled reactant forests has limit points. We use $(\,;\cdot)$ to denote the conditional law for the reactant with fixed realization of the catalyst environment.

PROPOSITION 4.1 (Existence of reactant limit forests). *Let $\tilde{\eta}^{\mathrm{tot},n}$ and $X$ satisfy (4.1). Then the family $\{(\tilde{\xi}^{\mathrm{for},n}; \tilde{\eta}^{\mathrm{tot},n}); n \in \mathbb{N}\}$ is relatively compact.*

To characterize the limit $(Y^{\mathrm{for}}; X)$ of forests $\{(\tilde{\xi}^{\mathrm{for},n}; \tilde{\eta}^{\mathrm{tot},n}), n \in \mathbb{N}\}$ we use the contour process and genealogical point process representations.

4.1.1. *The reactant limit contour in a fixed catalyst medium.* In order to obtain convergence of the rescaled reactant contour processes to a diffusion process we restrict our consideration only up to the first time the catalyst mass falls below an arbitrarily small threshold $\delta > 0$. Let

$$(4.2) \qquad \tilde{T}^{\delta,n} := \inf\{t \geq 0 : \tilde{\eta}_t^{\mathrm{tot},n} \leq \delta\}.$$

Recall that $Q_t$ defined in (2.2) takes a forest and cuts off the portion that lies above height $t$. Let $\tilde{\xi}^{\mathrm{for},\delta,n}$ denote the reactant trees of the $\tilde{\xi}^{\mathrm{for},n}$ cut off at $\tilde{T}^{\delta,n}$ (rather than at $\tilde{T}^{\delta,0}$)

$$\tilde{\xi}^{\mathrm{for},\delta,n} := Q_{\tilde{T}^{\delta,n}}(\tilde{\xi}^{\mathrm{for,n}})$$

and $\tilde{C}^{\delta,n}$ its contour process

$$\tilde{C}^{\delta,n} = (\tilde{C}_u^{\delta,n})_{u \geq 0} := \mathcal{C}(\tilde{\xi}^{\mathrm{for},\delta,n}; 2n).$$

The first result describes the behavior of the limit of the rescaled reactant forests in a catalytic environment that is stopped at $T^{\delta,n}$.

THEOREM 1 (Limit of the reactant contour process). *Assume (4.1), and fix $\delta > 0$. Consider the operator $(A^\delta, \mathcal{D}(A^\delta))$ with*

$$(4.3) \qquad A^\delta f(c) := 2\left(\frac{1}{X_c}f'\right)'(c)$$

*with reflection on the boundary, that is with domain*

$$(4.4) \qquad \mathcal{D}(A^\delta) := \left\{ f \in C^1([0,\tau^\delta]) : \frac{1}{X}. f' \in C^2([0,\tau^\delta]), f'|_{\{0,\tau^\delta\}} \equiv 0 \right\}.$$

*Then the following holds:*

(i) *The $(A^\delta, \mathcal{D}(A^\delta))$-martingale problem is well posed.*
(ii) *If $\zeta^\delta$ is the solution of the $(A^\delta, \mathcal{D}(A^\delta))$-martingale problem then*

$$(4.5) \qquad (\tilde{C}^{\delta,n}; \tilde{\eta}^{tot,n}) \underset{n\to\infty}{\Longrightarrow} (\zeta^\delta; X).$$

REMARK 4.1. By standard arguments (see, e.g., [25]) an operator in divergence form such as $A^\delta$ specifies a unique diffusion $\zeta^\delta$ with a well defined local-time process $\ell^t(\zeta^\delta) := (\ell_u^t(\zeta^\delta))_{u\geq 0}$ at all levels $t \geq 0$.

Recall that $\mathcal{T} : C^0_{\mathbb{R}^+}[0,\infty) \to \mathbb{T}^{root,lin}$ maps an excursion to a rooted linearly ordered compact $\mathbb{R}$-tree. Together with Theorem 1 we can immediately conclude the following.

COROLLARY 4.2. *For all $\delta > 0$, and limit points $Y^{for}$ of $\{(\tilde{\xi}^{for,n}; \tilde{\eta}^{tot,n}); n \in \mathbb{N}\}$*

$$(4.6) \qquad (Q_{\tau^\delta}(Y^{for}); X) \overset{d}{=} \mathcal{T}(\zeta^\delta; X).$$

This allows us to uniquely identify the reactant limit forest.

THEOREM 2 (The reactant limit forest exists). *There exists $Y^{for} \in \mathbb{T}^{root}$ such that*

$$(\tilde{\xi}^{for,n}; \tilde{\eta}^{tot,n}) \underset{n\to\infty}{\Longrightarrow} (Y^{for}; X).$$

This now raises the question whether we could let $\delta \to 0$ and obtain a limiting diffusion. Recall from (3.4) that $\tau^0$ is the first time the catalyst becomes extinct, and $\rho^0$ is the first time the reactant becomes extinct. As stated in (3.5), the probability of the reactant going extinct before the catalyst is positive. In case the reactant dies before the catalyst $\{\rho^0 < \tau^0\}$, the reactant is always in a strictly positive catalytic background. Hence, on this event, the limiting contour process is well defined by $X_t \in (0,\infty)$, for all $t \in [0, \rho^0]$, as a diffusion with generator

$$Af(c) := 2 \left( \frac{1}{X_c} f' \right)'(c), \qquad f \in \mathcal{D}(A),$$

on the domain

$$\mathcal{D}(A) := \left\{ f \in C^1([0,\tau^0]) : \frac{1}{X_\cdot} f' \in C^2_{[0,\tau^0]}[0,\infty), f'(0) = 0 \right\}.$$

However, on the event that the catalyst dies before the reactant $\{\rho^0 > \tau^0\}$, it is clear that there does not exist a well-defined limiting operator for $A^\delta$ as $\delta \to 0$. As will be shown in Corollary 4.7, on this event it is not possible to define a *diffusion* process $\zeta^0$ such that the tree below its path has the distribution of the limiting forest $Y^{for}$.

4.1.2. *The limit reactant point process in a fixed catalytic medium.* The point process representation of the genealogy is useful in those cases when we can give a simple description of its distribution. For both the catalyst and reactant genealogy one can show that given the total mass it can be represented by a Poisson point process with a specified intensity measure. Our next result describes the reactant limit genealogical point process.

THEOREM 3 (The reactant limit genealogical point process). *Assume* (4.1). *For any $t > 0$, and any sequence $(t_n > 0)_{n \in \mathbb{N}}$ such that $t_n \underset{n \to \infty}{\longrightarrow} t$,*

$$(\tilde{\Xi}^{t_n, n}; \tilde{\eta}^{\mathrm{tot}, n}) \underset{n \to \infty}{\Longrightarrow} (\pi^t; X),$$

*where for a given path $X = (X_t)_{t \geq 0}$, the point process $\pi^t$ consists of:*

*– rate $(\int_0^t X_s \, ds)^{-1}$ Poisson point process at height 0 whose points represent a separation of distinct trees in the forest,*

*– Poisson point process on $\mathbb{R}^+ \times \mathbb{R}^+$ whose intensity measure conditioned on $Y_t$ is*

$$(4.7) \qquad \aleph^{\zeta, t}[d\ell \times dh] = \mathbf{1}_{[0, Y_t]}(\ell) \, d\ell \otimes \mathbf{1}_{(0, t \wedge \tau^0 \wedge \rho^0)}(h) \frac{X_h}{(\int_h^t X_s \, ds)^2} \, dh.$$

We have only constructed a map $\mathcal{P}^t(\cdot; \varsigma)$ which sends an excursion of finitely many local maxima to a genealogical point process. To extend this map to excursions of infinitely many local maxima one would either proceed via a limit procedure or look at the point process of maximal depths of downward excursions. We next relate the reactant limit genealogical point process to the point process of depths of downward excursions corresponding to the limit contour.

PROPOSITION 4.3 (Relation between the limit point process and $\zeta^\delta$). *Fix $\delta > 0$, and let $\zeta^\delta$ be the solution of the well-posed martingale problem $(A^\delta, \mathcal{D}(A^\delta))$. Given $X := (X_s)_{s \geq 0}$, $\delta > 0$ and $t \in (0, \tau^\delta)$, the point process $\pi^t$ equals in distribution the point process of maximal depths of downward excursions of $\zeta^\delta$ below level $t$ indexed by the local-time process $(\ell_u^t(\zeta^\delta))_{u \geq 0}$ of $\zeta^\delta$.*

4.2. *Joint law of the catalyst and reactant.* The result of Theorem 2 can be extended to joint convergence of the catalyst and reactant rescaled forests, where for pairs of forests we use the product of Gromov–Hausdorff topologies.

THEOREM 4 [Limit of (catalyst, reactant)-forest]. *We have*

$$(4.8) \qquad (\tilde{\eta}^{\mathrm{for}, n}, \tilde{\xi}^{\mathrm{for}, n}) \underset{n \to \infty}{\Longrightarrow} (X^{\mathrm{for}}, Y^{\mathrm{for}}).$$

The result in Theorem 4 becomes relevant, for example, if we want to sample groups of individuals from both populations and to investigate the pair of matrices of mutual genealogical distances of sampled individuals.

4.3. *Differences between the reactant and "classical" forest.* Let $|\beta| = (|\beta|_u)_{u \geq 0}$ be a reflected standard branching Brownian motion, and for all $t \geq 0$, $\ell^t := (\ell^t_u)_{u \geq 0}$ its local time at level $t$. As before let $\ell^{-1}(s) := \inf\{u \geq 0 : \ell^0_u = s\}$. In the following we refer to

$$Z^{\text{for}} := \mathcal{T}(2|\beta|_{\cdot \wedge \ell^{-1}(1)})$$

as the *Brownian forest.* Recall from (3.6) together with (3.7) that $Z^{\text{for}}$ equals in distribution the catalyst forest $X^{\text{for}}$.

We now discuss the qualitative difference between the limiting forest $Y^{\text{for}}$ (with stochastically evolving branching rates), and the Brownian forest $Z^{\text{for}}$ (i.e., with fixed branching rates). We shall focus on the following aspects:

- differences in the tree structure due to inhomogeneity of the random environment (evolving branching rates);
- behavior of the forest $Y^{\text{for}}$ at time $\tau^0$ of extinction of the random environment (given the reactant survives the catalyst).

These two items can be considered from *quenched* and *annealed* points of view. Of course, for a statistician observing populations the annealed point of view is more relevant.

4.3.1. *Differences due to the inhomogeneous environment.* We first describe how in the quenched picture the reactant forest $Y^{\text{for}}$ can be obtained from the standard Brownian forest $Z^{\text{for}}$ by height dependent "stretching" of the tree metric expressed in terms of the total mass process of the fixed catalyst medium $X$. The next result compares the distribution of genealogical distances using the point-process representation. Recall from (3.4) that $\tau^0$ is the first time that $Z^{\text{for}}$ becomes extinct. The proposition compares in law the distances between points at a fixed level $t$ in a reactant forest $\partial Q_t(Y^{\text{for}})$, with the distances at another level $s_t(t)$ in a standard Brownian forest $\partial Q_{s_t(t)}(Z^{\text{for}})$.

PROPOSITION 4.4 (Stretching the tree metric). *Assume we are given a Brownian forest* $(Z^{\text{for}}, d_{Z^{\text{for}}})$. *For any* $\tau^0 \in (0, \infty)$ *and continuous function* $x : [0, \tau^0] \to \mathbb{R}_+$, *for a fixed* $t \in (0, \tau^0)$ *define an increasing function* $s^x_t : [0, t] \to [0, \int_0^t x_s \, ds]$ *by*

$$s^x_t(h) := \int_{t-h}^t x_s \, ds$$

*and let* $(s_t^x)^{-1} \colon [0, \int_0^t x_s \, ds] \to [0, t]$ *be its inverse. Then the following holds:*

$$((\partial Q_t(Y^{\mathrm{for}}), d_{Y^{\mathrm{for}}}); X = x) \stackrel{d}{=} (\partial Q_{s_t^x(t)}(Z^{\mathrm{for}}), 2 \cdot (s_t^x)^{-1}(\tfrac{1}{2} d_{Z^{\mathrm{for}}})).$$

Theorem 3 allows us to compare certain statistics of the reactant limit forest with those of a "classical" Brownian forest. The latter appears as the scaling limit of the uniform discrete forest with a fixed number of vertices as the number of vertices tends to infinity and may therefore be seen as the uniformly distributed compact rooted forest.

The following result states that provided we compare the reactant limit forest with the "classical" Brownian forest which has the same expected number of trees, the reactant limit forest has a smaller probability that two "randomly" picked individuals at a fixed time belong to different family trees than in the "classical" forest.

What does it mean to "pick individuals at random?" Recall that the map $\mathcal{T}$ takes an excursion and sends it to a rooted compact $\mathbb{R}$-tree. Fix $(T, \rho) \in \mathbb{T}^{\mathrm{root}}$, and let $h(T, \rho)$ denote the height of $T$ in terms of the largest distance of points in $T$ from the root. Assume that there exists an excursion $\zeta \in C^0_{\mathbb{R}^+}[0, \infty)$ with $\mathcal{T}(\zeta) = (T, \rho)$ such that the local-time processes $\ell^t(\zeta) = (\ell^t(\zeta)_s)_{s \geq 0}$ exist for all $t \in (0, h(T, \rho))$, and denote the inverse by $(\ell^t(\zeta))^{-1}$, that is, $(\ell^t(\zeta))^{-1}(x) := \inf\{u \geq 0 \colon \ell^t_u(\zeta) \geq x\}$. Let $L^t(\zeta) := \lim_{s \to \infty} \ell^t_s(\zeta) < \infty$. Define $H^t \colon [0, 1] \to T$ and $\mu^t$ by

$$H^t(u) := (\ell^t(\zeta))^{-1}(u \cdot L^t(\zeta)), \qquad \mu^t := (H^t)_* \lambda,$$

where $\lambda$ is Lebesgue measure on $[0, 1]$ and $H^t_* \lambda$ denotes the push forward of the Lebesgue measure under $H^t$. Obviously, $\mu^t$ can be thought of as the "uniform" distribution on the set $\partial Q_t(T, \rho)$.

PROPOSITION 4.5 (Comparison of the reactant with the classical tree). *Let* $(X^{\mathrm{for}}, Y^{\mathrm{for}})$ *be the catalytic branching forests with branching rates* $b_1 = b_2 = 1$, *and let* $Z^{\mathrm{for}}$ *be the forest associated with a standard reflected Brownian motion* $2|\beta|$ *stopped at* $\ell^{-1}(z) = \inf\{u \geq 0 \colon \ell_u(2|\beta|) = z\}$ *for some* $z \geq 0$. *Fix* $t \in (0, h(Y^{\mathrm{for}}) \wedge h(Z^{\mathrm{for}}))$, *and let* $Y_t = Z_t$ *and* $\mu^{t,Y}$ *and* $\mu^{t,Z}$ *be the "uniform" distributions on* $\partial Q_t(Y^{\mathrm{for}})$ *and* $\partial Q_t(Z^{\mathrm{for}})$, *respectively. Assume further that the expected number of trees alive at height* $t$ *are the same in* $Y^{\mathrm{for}}$ *and* $Z^{\mathrm{for}}$, *i.e. that* $z(t)^{-1} = \mathbb{E}[(\int_0^t X_s \, ds)^{-1}]$. *Then*

$$\mathbb{E}\left[\int_{(\partial Q_t(Y^{\mathrm{for}}))^2} (\mu^{t,Y})^{\otimes 2}(du, du') \mathbf{1}\{d_{Y^{\mathrm{for}}}(u, u') = 2t\}\right]$$

$$\leq \mathbb{E}\left[\int_{(\partial Q_t(Z^{\mathrm{for}}))^2} (\mu^{t,Z})^{\otimes 2}(du, du') \mathbf{1}\{d_{Z^{\mathrm{for}}}(u, u') = 2t\}\right].$$

4.3.2. *On top of the reactant limit tree.* In the region where the catalytic environment is about to go extinct the reactant has very few branching events. Hence the structure of the forest cannot be compared to the Brownian forest by a finite map of distance stretching. As a consequence the reactant limit forest will have infinite $\mathbf{l}^2$-length which is in contrast to the Brownian forest which is known to have finite $\mathbf{l}^2$-length.

To be more precise, let for a given $(T, d) \in \mathbb{T}$ and a finite subset $T' \subseteq T$,

$$\widetilde{\mathbf{l}}^2((T, d), T') := \sum_{a, b \in T' \,:\, ]a, b[ \cap T' = \varnothing} (d(a, b))^2.$$

We then say that $(T, d)$ is of finite (infinite) $\mathbf{l}^2$-*length* if there exists a finite (infinite) random variable $\mathbf{l}^2((T, d))$ such that for every sequence $(T'_n)_{n \in \mathbb{N}}$ of finite subsets of $T$ with $\inf\{\varepsilon > 0 : (T'_n)^\varepsilon = T\} \underset{n \to \infty}{\longrightarrow} 0$,

$$(4.9) \qquad \widetilde{\mathbf{l}}^2((T, d), T'_n) \underset{n \to \infty}{\longrightarrow} \mathbf{l}^2((T, d)) \qquad \text{in probability,}$$

where $A^\varepsilon$ denotes the $\varepsilon$-neighborhood of $A \subseteq T$.

LEMMA 4.6 (Infinite versus finite $\mathbf{l}^2$-length). *Let* $(X^{\mathrm{for}}, Y^{\mathrm{for}})$ *be the catalytic forests with branching rates* $b_1 = b_2 = 1$, *and* $Z^{\mathrm{for}}$ *a standard Brownian forest. Then*

$$\mathbb{P}\{\mathbf{l}^2(Y^{\mathrm{for}}) = \infty\} = \mathbb{P}\{\tau^0 \leq \rho^0\} > 0,$$

*while*

$$\mathbf{l}^2(Z^{\mathrm{for}}) < \infty, \qquad a.s.$$

From Lemma 4.6 we immediately conclude that there cannot be a diffusion $\zeta^0$ such that the forest below the path of a realization of $\zeta^0$ equals in distribution the reactant limit forest.

COROLLARY 4.7. *Let* $Y^{\mathrm{for}}$ *be the reactant limit forest, and* $\zeta$ *be a random element in* $C^0_{\mathbb{R}+}[0, \infty)$ *with* $\sup\{t \geq 0 : \zeta(t) \neq 0\} < \infty$, *a.s., such that* $\mathcal{T}(\zeta)$ *equals in distribution* $Y^{\mathrm{for}}$. *Then given that* $\{\tau^0 \leq \rho^0\}$ *the total quadratic variation of* $\zeta$ *is infinite. In particular, on the event* $\{\tau_0 < \rho_0\}$, $Y^{\mathrm{for}}$ *cannot be associated with the path of a diffusion process.*

4.4. *Extensions to more general processes.* In this paper we are mainly interested in the catalytic branching model as one of the three principal types of neutral branching models, as pointed out in the Introduction. If one assumed a perspective of branching processes in a random medium one could place the emphasis differently and start with more general Markovian dynamics for the catalyst process. One would assume that the scaling limit of the total mass of the catalyst and its genealogical tree exist and that the limit processes of the total mass of both the catalyst and reactant have paths in $C_{[0,\infty)}([0,\infty))$ with compact and connected support. Under these assumptions the techniques we developed in our proofs would imply convergence of the quenched total mass process of the reactant to a limit diffusion process, convergence of its quenched genealogical forest to a limit forest, as well as the convergence of the annealed bivariate objects. Similarly, Proposition 4.4 on stretching of the tree metric also remains valid.

Interesting examples could be provided by catalysts which are total mass dependent branching processes, that is time inhomogeneous models with branching rate $h(\eta_t^{\mathrm{tot}})$ where $h[0,\infty) \to [0,\infty)$ is a locally Lipschitz function with $h(0) = 0$.

If we are willing to give up results on the bivariate dynamics, we could also allow for the following. Suppose we are given the dynamics of a process of masses converging to a continuous path limit and still have all quenched results for the reactant forest. That is, assume that we have some scaling function $r : \mathbb{N} \to \mathbb{R}^+$ such that for all $T > 0$

$$\sup_{t \le T} |\eta_{r(n)t}^{\mathrm{tot},n} - X_t| \underset{n \to \infty}{\longrightarrow} 0$$

holds for some $X \in C_{[0,\infty)}([0,\infty))$ with $X^{-1}\{(0,\infty)\}$ almost surely a finite interval. Then Theorems 1, 2 and 3 hold again as well as Proposition 4.4. More generally the branching rate of the reactant could be $g(\eta_{r(n)}^{\mathrm{tot},n})$ for continuous functions $f$ vanishing at most at zero.

## 5. Proof of tightness (Proposition 4.1).

Throughout this section we assume that realizations of $\tilde{\eta}^{\mathrm{tot},n}, n \in \mathbb{N}$, and of $X$ on a single probability space are fixed such that (4.1) holds.

To prepare the proof of Proposition 4.1 we specify a compact set of rooted compact $\mathbb{R}$-trees which we will use for proving tightness. It will be helpful here to use genealogical terminology. For each $(T, \rho) \in \mathbb{T}^{\mathrm{root}}$ denote by $\partial T := T \setminus \bigcup_{x \in T} [\rho, x[$ the set of *leaves* of $T$, and by $h(T, \rho) := \sup_{x \in T} d(\rho, x)$ the *height* of the tree. We introduce the set $A_{t-\varepsilon}^t(T, \rho)$ of *ancestors* of the time $t$ population a time $\varepsilon$ back as follows: for all $\varepsilon > 0$ let

$$S_\varepsilon(T, \rho) := \{\rho\} \cup \{x \in T : \exists y \in T, \text{ with } x \in [\rho, y], \text{ and } d(x, y) \ge \varepsilon\}$$

be an $\varepsilon$-trimming of the tree $(T, \rho)$. In particular, if $h(T, \rho) \le \varepsilon$ then $S_\varepsilon(T, \rho) = \{\rho\}$. Then for given $(T, \rho) \in \mathbb{T}^{\mathrm{root}}$ let

$$A_{t-\varepsilon}^t(T, \rho) := S_\varepsilon(Q_t(T, \rho)) \cap \partial Q_{t-\varepsilon}(T, \rho)$$

be the set of ancestors at time $t - \varepsilon$ of those individuals that are alive at time $t$. It follows from Lemma 2.6 in [13] that $S_\varepsilon(T, \rho) \in \mathbb{T}^{\mathrm{root}}_{\mathrm{fin}}$ and hence $\#A_{t-\varepsilon}^t(T, \rho) < \infty$.

The following result says that the family of rooted compact $\mathbb{R}$-trees whose number of ancestors a positive time $\varepsilon$ back are uniformly bounded is precompact.

LEMMA 5.1 (A compact set of rooted $\mathbb{R}$-trees). *Fix a family* $\{L^m; m \in \mathbb{N}\}$ *of positive integers. Then the set*

$$
(5.1) \quad \Gamma := \Big\{ (T, \rho) \in \mathbb{T}^{\mathrm{root}};
$$
$$
\sum_{1 \le k \le \lfloor 2^{(m+1)} h(T,\rho) \rfloor} \# A_{(k-1)/2^{m+1}}^{k/2^{m+1}}(T, \rho) \le L^m, \forall m \in \mathbb{N} \Big\}
$$

*is compact in* $\mathbb{T}^{\mathrm{root}}$.

PROOF. Recall from Lemma 2.5 in [13] that a family $\Gamma$ is pre-compact in $(\mathbb{T}^{\mathrm{root}}, d_{\mathrm{GH}^{\mathrm{root}}})$ if for every $\varepsilon > 0$ there exists a positive integer $n(\varepsilon) < \infty$ such that each $T \in \Gamma$ has an $\varepsilon$-net with at most $n(\varepsilon)$ points.

Fix $\varepsilon > 0$ and let $(T, \rho) \in \Gamma$. Then with $m_0 := m_0(\varepsilon) = \lfloor -\log_2(\varepsilon) \rfloor$,

$$
N(\varepsilon) := \bigcup_{k \in \mathbb{N}} A_{(k-1)2^{-(m_0+1)}}^{k2^{-(m_0+1)}}(T, \rho)
$$

is an $\varepsilon$-net. Indeed, if $x \in T$, then for $B_\varepsilon(x) := \{y \in T : d(x, y) < \varepsilon\}$

$$
x \in B_\varepsilon(A_{2^{-(m_0+1)} \lfloor 2^{(m_0+1)} h(x)-1 \rfloor}^{2^{-(m_0+1)} \lfloor 2^{(m_0+1)} h(x) \rfloor}).
$$

Moreover, since $T \in \Gamma$, $\# N(\varepsilon) \le L^{m_0} = n(\varepsilon)$. Since $\Gamma$ is closed, compactness follows. □

PROOF OF PROPOSITION 4.1. We need to show that for all $\gamma > 0$, we can find a compact set $\Gamma := \Gamma_\gamma \subset \mathbb{T}^{\mathrm{root}}$ such that for all $n \in \mathbb{N}$,

$$
(5.2) \quad \mathbb{P}(\{\tilde{\xi}^{\mathrm{for},n} \in \Gamma\} | \tilde{\eta}^{\mathrm{tot},n}) \ge 1 - \gamma,
$$

provided that $\tilde{\eta}^{\mathrm{tot},n}$ and $X$ satisfy (4.1).

As candidates for such a compact set we consider sets of the form (5.1). To verify (5.2) calculate for fixed nondecreasing positive integers $\{L^m; m \in \mathbb{N}\}$,

$$
\mathbb{P}\Big\{ \sum_{1 \le k \le \lfloor 2^{(m+1)} \tilde{T}^{0,n} \rfloor} \# A_{2^{-(m+1)}(k-1)}^{2^{-(m+1)}k}(\tilde{\xi}^{\mathrm{for},n}) \le L^m, \forall m \in \mathbb{N} \Big| \tilde{\eta}^{\mathrm{tot},n} \Big\}
$$

$$
\ge 1 - \sum_{m \ge 1} \mathbb{P}\Big\{ \sum_{1 \le k \le \lfloor 2^{(m+1)} \tilde{T}^{0,n} \rfloor} \# A_{(k-1)2^{-(m+1)}}^{k2^{-(m+1)}}(\tilde{\xi}^{\mathrm{for},n}) \ge L^m \Big| \tilde{\eta}^{\mathrm{tot},n} \Big\}
$$

$$
\ge 1 - \sum_{m \ge 1} \frac{1}{L^m} \mathbb{E}\Big[ \sum_{1 \le k \le \lfloor 2^{(m+1)} \tilde{T}^{0,n} \rfloor} \# A_{(k-1)2^{-(m+1)}}^{k2^{-(m+1)}}(\tilde{\xi}^{\mathrm{for},n}) \Big| \tilde{\eta}^{\mathrm{tot},n} \Big].
$$

By Theorem 5.1 in [7], $(\#A_{t-\varepsilon}^t; \tilde{\eta}^{\mathrm{tot},n})_{\varepsilon \geq 0}$ is the number of blocks in a time-inhomogeneous Kingman coalescent at rate $1/\tilde{\eta}_{t-\varepsilon}^{\mathrm{tot},n}$. Moreover, if $K := (K_s)_{s \geq 0}$ is the number of blocks in a Kingman coalescent then $\lim_{\varepsilon \downarrow 0} \varepsilon K_\varepsilon = 2$ almost surely as well as $\lim_{\varepsilon \downarrow 0} \varepsilon \mathbb{E}(K_\varepsilon) = 2$, see page 27 of [4]. Hence there exists a constant $c > 0$ such that

$$
\begin{aligned}
(5.3) \quad & \mathbb{P}\left\{ \sum_{1 \leq k \leq \lfloor 2^{(m+1)} \tilde{T}^{0,n} \rfloor} \# A_{2^{-(m+1)}(k-1)}^{2^{-(m+1)}k}(\tilde{\xi}^{\mathrm{for},n}) \leq L^m, \forall m \in \mathbb{N} \Big| \tilde{\eta}^{\mathrm{tot},n} \right\} \\
& \geq 1 - c \sum_{m \geq 1} \frac{1}{L^m} \sum_{1 \leq k \leq \lfloor 2^{(m+1)} \tilde{T}^{0,n} \rfloor} \frac{1}{\int_{(k-1)2^{-(m+1)}}^{k2^{-(m+1)}} \tilde{\eta}_s^{\mathrm{tot},n} \, ds} \\
& \geq 1 - c \sum_{m \geq 1} \frac{1}{L^m} 2^{(m+2)} \Big( \min_{s \in [0, \tilde{T}^{0,n} - 2^{-(m+2)}]} \tilde{\eta}_s^{\mathrm{tot},n} \Big)^{-1} \lfloor 2^{(m+1)} \tilde{T}^{0,n} \rfloor \\
& \xrightarrow[n \to \infty]{} 1 - c \sum_{m \geq 1} \frac{1}{L^m} 2^{(m+2)} \Big( \min_{s \in [0, \tau^0 - 2^{-(m+2)}]} X_s \Big)^{-1} \lfloor 2^{(m+1)} \tau^0 \rfloor.
\end{aligned}
$$

Here we have used that for all $m \in \mathbb{N}$ and $k \in \{1, \ldots, \lfloor 2^{(m+1)} \tilde{T}^{0,n} \rfloor\}$,

$$
\begin{aligned}
\int_{(k-1)2^{-(m+1)}}^{k2^{-(m+1)}} \tilde{\eta}_s^{\mathrm{tot},n} \, ds & \geq \int_{(k-1)2^{-(m+1)}}^{k2^{-(m+1)} - 2^{-(m+2)}} \tilde{\eta}_s^{\mathrm{tot},n} \, ds \\
& \geq 2^{-(m+2)} \cdot \min_{s \in [0, \tilde{T}^{0,n} - 2^{-(m+2)}]} \tilde{\eta}_s^{\mathrm{tot},n}.
\end{aligned}
$$

Since $\tau^0 < \infty$ and $\min_{s \in [0, \tau^0 - 2^{-(m+2)}]} X_s > 0$, for all $m \in \mathbb{N}$ and all realizations of $X$, we can choose the family of positive integers $\{L^m, m \in \mathbb{N}\}$ suitably large, for example, $L^m > (4/q)^m (\min_{s \in [0, \tau^0 - 2^{-(m+2)}]} X_s)^{-1}$ with $0 < q = \gamma/(\gamma + 1) < 1/2$, such that the right-hand side of (5.3) is greater than or equal to $1 - \gamma$, and (5.2) follows. $\square$

## 6. Proof of the limit for the reactant contour (Theorem 1).

In this section we give the proof of Theorem 1. We first show that given a realization of the catalyst total mass process the reactant contour process is associated with a Markov process, and we derive its generator. We also give a representation of this Markov process as a *random evolution* process. By this we mean that it moves at constant velocity for a random time, then changes the sign of its velocity and proceeds at constant velocity for a random time again. We use this representation to prove that the suitably rescaled contour processes of the truncated reactant forest converge toward a limit contour process which is characterized as the solution of a well-posed martingale problem.

Throughout this section a realization $\eta^{\text{tot}} := (\eta_s^{\text{tot}})_{s \geq 0} \in D_{\mathbb{N}}[0, \infty)$ of the catalyst path is fixed. Recall the reactant contour process $C := (C_u)_{u \geq 0}$ from Definition 3.3, and let its *slope process* $V := (V_u)_{u \geq 0}$ be defined by

$$V_u := \tfrac{1}{2} \operatorname{slope}(C_u) \in E_{\text{slope}}$$

with

$$E_{\text{slope}} := \{-1, +1\}, \qquad E_{\text{cont}} := [0, T^{0,1}],$$

where the slope of the root, branch points and leaves are defined in such a way that $(V_u)_{u \geq 0}$ has cádlág paths and $T^{0,1} := \inf\{t \geq 0 \colon \eta_t^{\text{tot}} = 0\}$.

The next result states that the functional pairing of the height of the contour with its slope is a Markov process.

LEMMA 6.1 (Markov property of the contour process). *The process* $(C, V) := (C_u, V_u)_{u \geq 0}$ *on* $E_{\text{cont}} \times E_{\text{slope}}$ *is a Markov process whose generator is the closure of the operator*

$$(6.1) \qquad Af(c, v) = 2v \frac{\partial}{\partial c} f(c, v) + \eta_c^{\text{tot}} [f(c, -v) - f(c, v)]$$

*for all* $f \in \mathcal{D}(A)$, *where*

$$(6.2) \qquad \mathcal{D}(A) = \left\{ f \in C^{1,0}(E_{\text{cont}} \times E_{\text{slope}}) \colon \frac{\partial f}{\partial c} \Big|_{\partial(E_{\text{cont}}) \times E_{\text{slope}}} \equiv 0 \right\}.$$

PROOF. Recall that $\mathcal{C}(\cdot; \sigma) \colon \mathbb{T}_{\text{fin}}^{\text{root,lin}} \to C_{\mathbb{R}+}[0, \infty)$ maps a rooted $\mathbb{R}$-tree to an excursion from Figure 2. We first show that the lengths of the line segments of the contour process are distributed as a sequence of independent variables stopped at the first time their alternating sum falls below 0, and then we use this to obtain the Markov property and to identify the generator.

STEP 1 (Sequence of lengths of the line segments). Recall from Figure 2 the assignment of a contour process to the representative of a Galton–Watson family forest embedded in the plane. Each piece of the contour process with constant slope sign corresponds to a sum of a number of lifetimes, so the distributional relationship between the different line segments is not obvious. We shall make use of the fact that the reactant process can be represented in two ways without changing the distribution of its total mass process or the genealogical distances between individuals, either as continuous-time *binary Galton–Watson process*, or as a *birth-and-death process*:

- In the continuous-time Galton–Watson process the branch points occur at rate $\eta_t^{\text{tot}}$, and the number of offspring at each branch point is 0 or 2 with equal probability.

- In the birth-and-death process each individual dies at rate $\frac{1}{2}\eta_t^{\mathrm{tot}}$, and during its lifetime gives birth to new offspring at rate $\frac{1}{2}\eta_t^{\mathrm{tot}}$.

If in the Galton–Watson process at each birth time we choose to identify the life of one of the offspring as a continuation of the life of its parent, we obtain the birth-and-death process.

The family forest of the Galton–Watson process has a canonical planar embedding coming from a linear order on its vertices induced by the linear order of the tree as in Figure 3(a), while the family forest of a birth-and-death process with that same linear order on the vertices has two canonical planar embeddings. In Figure 3(b) we always choose to identify the continuation of the life of the parent with the life of the offspring of higher linear order, and the branch of the offspring is always drawn to the left of the branch of the parent. In Figure 3(c) we always choose to identify the continuation of the life of the parent with the life of the offspring of lower linear order, and the branch of the offspring is always drawn to the right of the branch of the parent.

The key observation now is that since all three planar embeddings respect the same linear order on the vertices they also have the same contour process (compare Figure 2).

In the birth-and-death process the line segments of constant slope correspond to a lifetime of exactly one individual. With the parent identification as in Figure 3(b) each line segments of negative slope corresponds to a lifetime of an individual, while in the parent identification as in Figure 3(c) each line segment of positive slope corresponds to a lifetime of an individual. Since lifetimes of individuals are independent this implies the constant slope line segments in the contour process are distributed as a sequence of independent variables stopped when their alternating sum first falls below 0.

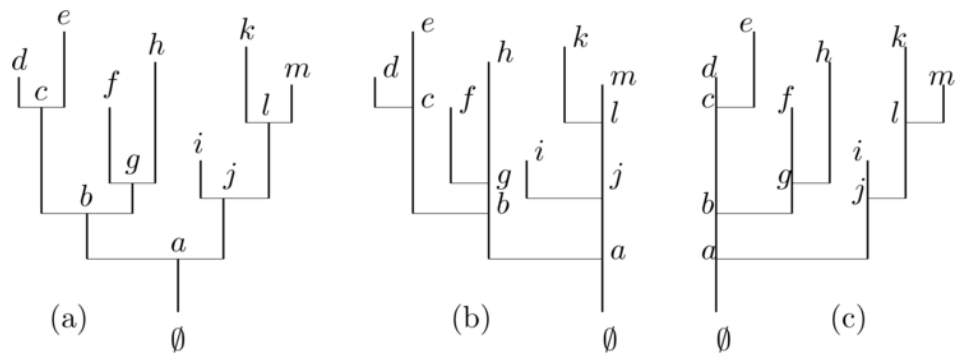

FIG. 3. *Illustrates three different planar embeddings of the reactant family forest with the same linear order and a common contour process.*

STEP 2 (Identification of the generator). Given the branch rates $(\eta_t^{\text{tot}})_{t\geq0}$, the law of the length $l$ of a lifetime of an individual is

$$\mathbb{P}\{l > t\} = e^{-\int_{t^*}^{t^*+t} \eta_s^{\text{tot}}\, ds},$$

where $t^*$ is the time of birth of that individual. Hence the law of the length of each segment is the law of the first point of a Poisson process with rate $(\eta_t^{\text{tot}})_{t\geq t^*}$ or $(\eta_t^{\text{tot}})_{t\leq t^*}$ if the slope of the line segment is $+2$ and $-2$, respectively. The alternating sign of $V_u$ changes at rate $\eta_{C_u}^{\text{tot}}$ where $C_u$ is the current value of the contour process. In between the jumps of $V_u$, $C_u$ moves with speed of 2 units in the direction determined by $V_u$. Hence, the contour process paired with its slope is Markovian, and furthermore its generator agrees on $\mathcal{D}(A)$ with the operator given in (6.1) and (6.2). Since $\mathcal{D}$ is dense in $C(E_{\text{cont}} \times E_{\text{slope}})$ the generator is the closure of the operator $(A, \mathcal{D})$.  □

We next notice that $(C, V)$ is a *random evolution*, that is, a Markov process moving at constant speed in a direction whose sign changes stochastically. Specifically, for the pair $(C, V)$ the change in speed is a counting process whose rate is governed by the catalyst mass process $\eta^{\text{tot}}$.

LEMMA 6.2 (Random evolution representation). *Let $N := (N(u))_{u\geq0}$ be a unit rate Poisson process, and consider the following system*

$$(6.3) \qquad C_u = \int_0^u 2V_v\, dv, \qquad V_u = (-1)^{N\left(\int_0^u \eta_{C_v}^{\text{tot}}\, dv\right)}$$

*in* $\text{int}(E_{\text{cont}}) \times E_{\text{slope}}$ *and with reflection on* $\partial(E_{\text{cont}}) \times E_{\text{slope}}$. *There exists a unique random evolution satisfying this system, and its distribution is the same as the distribution of the contour process and its slope* $(C, V)$.

PROOF. See Chapter 12 of [12] for the definition of a random evolution. Existence and uniqueness of a random evolution satisfying the system (6.3) follow from standard theorems on existence and solutions of a system of stochastic differential equations with continuous local martingales as differentials (see Theorem 3.15 in [23]). Equality in distribution of this random evolution with the contour process and its slope follows by simply comparing the generator of this system to the generator obtained in the previous lemma.  □

We will also need the following easy consequence of the above arguments. Fix $\delta > 0$, and recall that $T^{\delta,1}$ from (4.2) is the first time that the catalyst process started from 1 individual drops below $\delta$. Also recall that $Q_{T^{\delta,1}}$ from (2.2) is the map that takes a tree and cuts the portion of its branches that lie above height $T^{\delta,1}$.

COROLLARY 6.3 (Truncated process). *If the contour $\mathcal{C}(\xi^{\text{for}}; 1)$ solves (6.3) then $\mathcal{C}(Q_{T^{\delta,1}}(\xi^{\text{for}}); 1)$ solves (6.3) but with the state space $E_{\text{cont}}$ replaced by $E_{\text{cont}}^{\delta} = [0, T^{\delta,1}]$.*

We next show that the rescaled random evolutions converge to the solution of the well-posed martingale problem stated in Theorem 1. We rely on averaging techniques established for random evolutions. Throughout we fix realizations of $\tilde{\eta}^{\text{tot},n}$, $n \in \mathbb{N}$ and $X$ on a single probability space such that (4.1) holds and we choose a truncation parameter $\delta > 0$.

STEP 0 (The rescaled random evolution system). Recall the rescaled reactant contour process $\tilde{C}^{\delta,n}$ from Definition 3.7. Let $\tilde{E}_{\text{cont}}^{\delta,n} := [0, \tilde{T}^{\delta,n}]$, and define its rescaled $\frac{1}{2}$-slope process $\tilde{V}^{\delta,n} := \frac{1}{2}\,\text{sign}(\text{slope}(\tilde{C}_{\cdot}^{\delta,n}))$.

Then Lemma 6.1 applied to the rescaled reactant populations implies that $(\tilde{C}^{\delta,n}, \tilde{V}^{\delta,n})$ is a Markov process whose generator is the closure of the operator

$$(6.4) \qquad \tilde{A}^{\delta,n} f(c,v) = 2nv \frac{\partial}{\partial c} f(c,v) + n^2 \tilde{\eta}_c^{\text{tot},n}[f(c,-v) - f(c,v)],$$

acting on all $f \in \mathcal{D}(\tilde{A}^{\delta,n})$, where

$$\mathcal{D}(\tilde{A}^{\delta,n}) = \left\{ f \in C^{1,0}(\tilde{E}_{\text{cont}}^{\delta,n} \times E_{\text{slope}}) : \frac{\partial f}{\partial c}\bigg|_{\{0, \tilde{T}^{\delta,n}\} \times E_{\text{slope}}} \equiv 0 \right\}.$$

Furthermore, the analogous argument to that of Lemma 6.2 implies that the distribution of the pair $(\tilde{C}^{\delta,n}, \tilde{V}^{\delta,n})$ has the same distribution as the random evolution which is the unique solution to the system

$$(6.5) \qquad \tilde{C}_u^{\delta,n} = n \int_0^u 2\tilde{V}_v^{\delta,n}\, dv, \qquad \tilde{V}_u^{\delta,n} = (-1)^{N\left(n^2 \int_0^u \tilde{\eta}_{\tilde{C}_v^{\delta,n}}^{\text{tot},n}\, dv\right)},$$

where $N$ is a unit rate Poisson process. In other words, the rescaled process $\tilde{C}^{\delta,n}$ evolves deterministically with speed $2n$, and changes the sign at rate $n^2$ times a counting process whose rate is governed by the rescaled catalyst mass process $\tilde{\eta}^{\text{tot},n}$.

STEP 1 (The velocity process). The convergence result stated in Theorem 1 relies on the fact that the velocity component $\tilde{V}^{\delta,n}$ evolves much faster than the contour component $\tilde{C}^{\delta,n}$, as it is clear from (6.5). Hence in the limit the velocity process will average out and can be replaced with its stationary measure.

If $\Gamma^n$ is the occupation-time measure of $\tilde{V}^{\delta,n}$ on $E_{\text{slope}}$, that is, for $u \geq 0$, and $v \in E_{\text{slope}}$

$$(6.6) \qquad \Gamma^n([0,u] \times \{v\}) := \int_0^u \mathbf{1}_{\{v\}}(\tilde{V}_{u'}^{\delta,n})\, du',$$

then it is clear from the description (6.5) of $\tilde{V}^{\delta,n}$ that

$$\Gamma^n \underset{n\to\infty}{\Longrightarrow} \Gamma = \lambda \otimes \pi,$$

where $\lambda$ denotes Lebesgue measure on $[0,\infty)$, and $\pi(\{1\}) = \pi(\{-1\}) = \frac{1}{2}$. In the limit the contour component will have spent on average half the time increasing and half the time decreasing.

STEP 2 (Averaging for martingale problems). The proof of Theorem 1 will rely on the following result taken from Theorem 2.1 in [17] and adapted to our specific situation in which the state spaces are compact.

PROPOSITION 6.4 (Stochastic averaging). *Suppose there is an operator* $A^\delta : \mathcal{D}(A^\delta) \subseteq C^0([0,\tau^\delta]) \to C^0([0,\tau^\delta] \times \{-1,1\})$ *such that:*

(i) *for all* $f \in \mathcal{D}(A^\delta)$ *there is a process* $\varepsilon^{\delta,f,n}$ *for which*

$$(6.7) \qquad \left( f(\tilde{C}_t^{\delta,n}) - \int_0^t (A^\delta f)(\tilde{C}_s^{\delta,n}, \tilde{V}_s^{\delta,n})\, ds + \varepsilon_t^{\delta,f,n} \right)_{t \geq 0}$$

*is a martingale,*

(ii) $\mathcal{D}(A^\delta)$ *is dense in* $C^0([0,\tau^\delta])$ *with respect to the uniform topology and*

(iii) *for* $f \in \mathcal{D}(A^\delta)$ *and* $T > 0$ *there exists* $p > 1$ *such that*

$$(6.8) \qquad \sup_n \mathbb{E}\left[ \int_0^T |A^\delta f(\tilde{C}_s^{\delta,n}, \tilde{V}_s^{\delta,n})|^p\, ds \right] < \infty$$

*and*

$$(6.9) \qquad \lim_{n\to\infty} \mathbb{E}\left[ \sup_{t \leq T} |\varepsilon_t^{\delta,f,n}| \right] = 0.$$

*Then for* $\Gamma^n$ *defined as in (6.6)*

$$\{(\tilde{C}^{\delta,n}, \Gamma^n)\}_{n\in\mathbb{N}} \text{ is relatively compact in } \mathcal{D}([0,\infty) \times \mathcal{M}([0,\infty) \times \{-1,1\})),$$

*where* $\mathcal{M}([0,\infty) \times \{-1,1\})$ *is the space of all measures on* $[0,\infty) \times \{-1,1\}$ *for which* $\mu([0,u) \times \{-1,1\}) = u$, *and for any limit point* $(\zeta^\delta, \pi)$ *there exists a filtration such that*

$$\left( f(\zeta_t^\delta) - \int_0^t \sum_{v\in\{-1,1\}} A^\delta f(\tilde{C}_s^{\delta,n}, v)\, ds\, \pi\{v\} \right)_{t \geq 0}$$

*is a martingale with respect to this filtration, for all* $f \in \mathcal{D}(A^\delta)$.

STEP 3 (The averaging-martingale problem). Recall the operator $(A^\delta, \mathcal{D}(A^\delta))$ from (4.3) and (4.4) which is the generator for the rescaled contour process $\tilde{C}^{\delta,n}$. Our goal is to show that all three assumptions (i)–(iii) of Proposition 6.4 above are satisfied.

(i) We first show that we can define small error processes $\{\varepsilon^{\delta,f,n}; f \in \mathcal{D}(A^\delta)\}$ such that

$$\left( f(\tilde{C}_t^{\delta,n}) - \int_0^t A^\delta f(\tilde{C}_s^{\delta,n}, v)\, ds + \varepsilon_t^{\delta,f,n} \right)_{t \geq 0}$$

is a martingale for all $f \in \mathcal{D}(A^\delta)$.

Notice that $f \in \mathcal{D}(A^\delta)$ if and only if there exists a function $g \in C^2([0,\tau^\delta])$ with $g(0) = g(\tau^\delta) = 0$ such that

(6.10) $$f(x) = f(0) + \int_0^x X_s g(s)\, ds$$

for all $x \geq 0$. Let $f \in \mathcal{D}(A^\delta)$ be of the form (6.10) and let

$$\tilde{f}^n(c) := f(0) + \int_0^c \tilde{\eta}_s^{\mathrm{tot},n} g\left( \frac{\tilde{T}^{\delta,n}}{\tau^\delta} s \right) ds.$$

Apply the operator $(\tilde{A}^{\delta,n}, \mathcal{D}(\tilde{A}^{\delta,n}))$ from (6.4) to functions $f^n$ given by

$$f^n(c, v) := \tilde{f}^n(c) + \frac{v}{n\tilde{\eta}_c^{\mathrm{tot},n}} (\tilde{f}^n)'(c)$$

to get

$$\tilde{A}^{\delta,n} f^n(c,v) = 2nv(\tilde{f}^n)'(c) - 2n^2 \tilde{\eta}_c^{\mathrm{tot},n} \frac{v}{n\tilde{\eta}_c^{\mathrm{tot},n}} (\tilde{f}^n)'(c) + 2v^2 \left( \frac{1}{\tilde{\eta}_c^{\mathrm{tot},n}} (\tilde{f}^n)' \right)'(c)$$

$$= 2\left( \frac{1}{\tilde{\eta}_c^{\mathrm{tot},n}} (\tilde{f}^n)' \right)'(c).$$

Let then for all $t \geq 0$,

$$\varepsilon_t^{\delta,f,n} := \tilde{f}^n(\tilde{C}_t^{\delta,n}) - f(\tilde{C}_t^{\delta,n}) + \frac{\tilde{V}_t^{\delta,n}}{n\tilde{\eta}_{\tilde{C}_t^{\delta,n}}^{\mathrm{tot},n}} (\tilde{f}^n)'(\tilde{C}_t^{\delta,n})$$

$$+ \int_0^t (A^\delta f - \tilde{A}^{\delta,n} f^n)(\tilde{C}_s^{\delta,n}, \tilde{V}_s^{\delta,n})\, ds.$$

Since $f^n(\tilde{C}_t^{\delta,n}, \tilde{V}_t^{\delta,n}) - \int_0^t \tilde{A}^{\delta,n} f^n(\tilde{C}_s^{\delta,n}, \tilde{V}_s^{\delta,n})\, ds$ is a martingale for all $f^n \in \mathcal{D}(\tilde{A}^{\delta,n})$, it follows that (6.7) holds for all $f \in \mathcal{D}(A^\delta)$.

(ii) We next show that the domain $\mathcal{D}(A^\delta)$ is dense in the space of continuous functions on $[0, \tau^\delta]$.

LEMMA 6.5 (Dense domain). *Fix $\delta > 0$ and $X \in C^0([0,\tau^\delta])$. Then the set of functions $\mathcal{F}$ defined by*

$$\mathcal{F} := \left\{ f : f(c) = C + \int_0^c X_{c'} g(c')\, dc'; C \in \mathbb{R}, g \in C^2([0,\tau^\delta]), g|_{\{0,\tau^\delta\}} \equiv 0 \right\}$$

*is dense in $C^0([0,\tau^\delta])$.*

PROOF.  It is well known that each continuous function on $[0, \tau^\delta]$ can be approximated by piecewise linear functions. It is therefore enough to show that any piecewise linear function can be approximated by functions in $\mathcal{F}$. This follows by continuity of $X$ and the fact that $X_u \geq \delta$, for all $u \in [0, \tau^\delta]$. $\square$

(iii) We finally verify the last point. It is standard to show that (6.8) holds for all $f \in \mathcal{D}(A^\delta)$, $T > 0$ and $p > 1$. Moreover, since $1/\tilde{\eta}^{\mathrm{tot},n}_{\tilde{C}^{\delta,n}_t}$ is bounded by $\frac{1}{\delta}$, for all $t \geq 0$, $\tilde{f}^n \to f$, and $\|\tilde{A}^{\delta,n} f^n - A^\delta f\| \leq |\tau^\delta - T^{\delta,n}| \|g'\| + |1 - \frac{T^{\delta,n}}{\tau^\delta}| \|g''\| \underset{n \to \infty}{\longrightarrow} 0$, (6.9) is satisfied as well.

Altogether we can apply Proposition 6.4 to the effect that the family of rescaled contours $\{\tilde{C}^{\delta,n}; n \in \mathbb{N}\}$ is relatively compact in law and any limit point satisfies the $(A^\delta, \mathcal{D}(A^\delta))$-martingale problem.

STEP 4 (Uniqueness of the limit-martingale problem).  We next show that the $(A^\delta, \mathcal{D}(A^\delta))$-martingale problem has a unique solution $\zeta^\delta$, for which we first need the following lemma which characterizes solutions of transforms of a reflecting Brownian motion. Recall that $X$ is a Feller diffusion started at $X_0 = 1$ and $\tau^\delta$ is the first time it falls below $\delta$. Let $s : [0, \tau^\delta] \to [0, \int_0^{\tau^\delta} X_u \, du]$ defined by

$$(6.11) \qquad\qquad s(x) := \int_0^x X_u \, du,$$

denote the *scale function* and denote its inverse by $s^{-1}$.

Consider the operator

$$(B^\delta, \mathcal{D}(B^\delta)) := (2X_{s^{-1}(\cdot)} f'', \{h \in C^2_{\mathbb{R}^+}[0, \infty) : h'|_{\{0, s(\tau^\delta)\}} \equiv 0\}).$$

LEMMA 6.6 (Relation to Brownian motion).  *Let $\delta > 0$. If $\zeta^\delta$ solves the $(A^\delta, \mathcal{D}(A^\delta))$-martingale problem, then*

$$B_t := s(\zeta^\delta_t), \qquad t \in [0, \infty),$$

*solves the $(B^\delta, \mathcal{D}(B^\delta))$-martingale problem.*

PROOF.  Fix $H \in C^2_{[0, s(\tau^\delta)]}[0, \infty)$ with $H'|_{\{0, s(\tau^\delta)\}} \equiv 0$. It is easy to check that then $H \circ s \in \mathcal{D}(A^\delta)$.

We therefore obtain that for all $H \in C^2_{\mathbb{R}^+}[0, \infty)$ such that $H'|_{\{0, s(\tau^\delta)\}} \equiv 0$ the process

$$\left( (H \circ s)(\zeta^\delta_t) - \int_0^t A^\delta(H \circ s)(\zeta^\delta_u) \, du \right)_{t \geq 0}$$

$$= \left( H(B_t) - \int_0^t 2X_{s^{-1}(B_u)} H''(B_u) \, du \right)_{t \geq 0}$$

is a martingale. Here we have used that $A^\delta(H \circ s)(\zeta) = 2(H' \circ s)'(\zeta) = 2X_{\zeta}.(H'' \circ s)(\zeta)$. □

PROPOSITION 6.7 (Uniqueness).  *The $(A^\delta, \mathcal{D}(A^\delta))$-martingale problem has a unique solution $\zeta^\delta$.*

PROOF OF PROPOSITION 6.7.  Assume that $\zeta^{\delta,1}$ and $\zeta^{\delta,2}$ are two solutions of the $(A^\delta, \mathcal{D}(A^\delta))$-martingale problem. By Theorem 5.6 in [26], since $(X_u)_{u \geq 0}$ is bounded away from zero, the $(B^\delta, \mathcal{D}(B^\delta))$-martingale problem has a unique solution. Hence

$$\int_{[0, \zeta^{\delta,1}]} X_u \, du \overset{d}{=} \int_{[0, \zeta^{\delta,2}]} X_u \, du.$$

Since the scale function is strictly increasing on $[0, \tau^\delta]$, the one-dimensional distributions of $\zeta^{\delta,1}$ and $\zeta^{\delta,2}$ agree. It follows from Theorem 4.4.2 in [12] that therefore $\zeta^{\delta,1} \overset{d}{=} \zeta^{\delta,2}$. □

STEP 5 (Conclusion).  We close the section by giving a proof of Theorem 1.

PROOF OF THEOREM 1.  We have shown in Step 2 that the sequence $(\tilde{C}^{\delta,n})_{n \in \mathbb{N}}$ is relatively compact, and that any limit point $\zeta^\delta$ of $\tilde{C}^{\delta,n}$ is a solution of the $(A^\delta, \mathcal{D}(A^\delta))$-martingale problem. This proves existence of a solution for all $\delta > 0$.

Furthermore, by Proposition 6.7, this martingale problem has only one solution. That is, the $(A^\delta, \mathcal{D}(A^\delta))$-martingale problem is well posed and if $\zeta^\delta$ is its unique solution (4.5) holds.

Finally, it follows from Theorem 4.4.2 in [12] that $\zeta^\delta$ has the Markov property. □

## 7. Proof of convergence of the rescaled reactant forests (Theorem 2).
In this section we verify that the sequence $\{\tilde{\xi}^{\mathrm{for},n}; n \in \mathbb{N}\}$ converges in law. So far we have shown in Section 5 that this sequence of laws is relatively compact. It therefore remains to show that there is a unique limit law. We first give a proof of Corollary 4.2.

PROOF OF COROLLARY 4.2.  Fix $\delta > 0$. By Proposition 4.1 there exists a random forest $Y^{\mathrm{for}}$ such that along a subsequence $(n_k)_{k \in \mathbb{N}}$,

$$(\tilde{\xi}^{\mathrm{for},n_k}; \tilde{\eta}^{n_k}) \underset{k \to \infty}{\Longrightarrow} (Y^{\mathrm{for}}; X)$$

and hence, by the joint continuity of the truncation map $Q.(\cdot) \colon \mathbb{R}^+ \times \mathbb{T}^{\mathrm{root}} \to \mathbb{T}^{\mathrm{root}}$, we have

$$(Q_{\tilde{T}^{\delta,n_k}}(\tilde{\xi}^{\mathrm{for},n_k}); \tilde{\eta}^{n_k}) \underset{k \to \infty}{\Longrightarrow} (Q_{\tau^\delta}(Y^{\mathrm{for}}); X).$$

On the other hand we have proven in Theorem 1 that

$$(\mathcal{C}(Q_{\tilde{T}^{\delta,n_k}}(\tilde{\xi}^{\mathrm{for},n_k}); 2n_k); \tilde{\eta}^{n_k}) =: (\tilde{C}^{\delta,n_k}; \tilde{\eta}^{n_k}) \underset{k\to\infty}{\Longrightarrow} (\zeta^\delta; X).$$

Finally, since the map $\mathcal{T}: C^0_{\mathbb{R}^+}[0,\infty) \to \mathbb{T}^{\mathrm{root}}$ is continuous, and for all $\xi^{\mathrm{for}} \in \mathbb{T}^{\mathrm{root}}_{\mathrm{fin}}$ we have $\mathcal{T}(C(\xi^{\mathrm{for}}; \sigma)) = \xi^{\mathrm{for}}$, we know that for all $\delta > 0$,

$$(Q_{\tilde{T}^{\delta,n_k}}(\tilde{\xi}^{\mathrm{for},n_k}); \tilde{\eta}^{n_k}) \underset{k\to\infty}{\Longrightarrow} (\mathcal{T}(\zeta^\delta); X),$$

which proves (4.6). □

PROOF OF THEOREM 2. It follows from Corollary 4.2 that for all limit forests $Y^{\mathrm{for}}$ of $\{\tilde{\xi}^{\mathrm{for},n}; n \in \mathbb{N}\}$, $(Q_{\tau^\delta}(Y^{\mathrm{for}}); X)$ equals in distribution $(\mathcal{T}(\zeta^\delta); X)$. We will use this to show that the family $\{\mathcal{L}[(Q_{\tau^\delta}(Y^{\mathrm{for}}); X)]; \delta > 0\}$ is a Cauchy sequence in the Prohorov metric of probability measures. For that fix $0 < \delta' < \delta$, and recall that the Prohorov metric satisfies (see, e.g., Theorem 3.1.2 in [12]),

$$d_{\mathrm{Pr}}(\mathcal{L}(\mathcal{T}(\zeta^\delta); X), \mathcal{L}(\mathcal{T}(\zeta^{\delta'}); X))$$
$$:= \inf_\mu \inf\{\varepsilon > 0 : \mu\{((T,\rho),(T',\rho')) : d_{\mathrm{GH}^{\mathrm{root}}}((T,\rho),(T',\rho')) \geq \varepsilon\} \leq \varepsilon\},$$

where the infimum is taken over all probability measures $\mu$ on $\mathbb{T}^{\mathrm{root}} \times \mathbb{T}^{\mathrm{root}}$ with marginals $\mathcal{L}(\mathcal{T}(\zeta^\delta); X)$ and $\mathcal{L}(\mathcal{T}(\zeta^{\delta'}); X)$.

A particular choice of $\mu$ is the following: pick an excursion $\zeta^{\delta'}$ according to $\mathcal{L}(\mathcal{T}(\zeta^{\delta'}))$, and recall that $\ell^{\tau_\delta}(\zeta^{\delta'})$ is the local-time process of $\zeta^{\delta'}$ at level $\tau_\delta$. Write $\zeta^{\delta', \downarrow \tau^\delta}$ for $\zeta^{\delta'}$ time changed by the inverse of $t \mapsto \int_0^t \mathbf{1}\{\zeta^{\delta'}_s \leq \tau^\delta\} ds$. That is, $\zeta^{\delta', \downarrow \tau^\delta}$ is $\zeta^{\delta'}$ with the sub-excursions above level $\tau^\delta$ excised and the gaps closed up. Then $\mathcal{T}(\zeta^{\delta', \downarrow \tau^\delta})$ equals in law $\mathcal{T}(\zeta^\delta)$ and clearly,

$$d_{\mathrm{GH}^{\mathrm{root}}}(\mathcal{T}(\zeta^{\delta'}), \mathcal{T}(\zeta^{\delta', \downarrow \tau^\delta})) = d_{\mathrm{H}}(\mathcal{T}(\zeta^{\delta'}), \mathcal{T}(\zeta^{\delta', \downarrow \tau^\delta})) = \tau^{\delta'} - \tau^\delta.$$

Hence,

$$d_{\mathrm{Pr}}(\mathcal{L}(\mathcal{T}(\zeta^\delta); X), \mathcal{L}(\mathcal{T}(\zeta^{\delta'}); X)) \leq \tau^{\delta'} - \tau^\delta,$$

which proves that $\{\mathcal{L}[(Q_{\tau^\delta}(Y); X)]; \delta > 0\}$ is Cauchy, as claimed.

Since $(\mathbb{T}^{\mathrm{root}}, d_{\mathrm{GH}^{\mathrm{root}}})$ is complete, the space of probability measures on $(\mathbb{T}^{\mathrm{root}}, d_{\mathrm{GH}^{\mathrm{root}}})$ is complete as well, and we can therefore find a $\mathbb{P} \in \mathcal{M}_1(\mathbb{T}^{\mathrm{root}})$ such that $\mathcal{L}(Q_{\tau^\delta}(T)) \underset{\delta \to 0}{\Longrightarrow} \mathbb{P}$, and Theorem 2 is proved. □

**8. Proofs for the reactant limit point process (Theorem 3).** In this section we derive the distribution of the reactant limit point process by taking a limit of the point processes for the finite reactant forests, as well as constructing the limit point process from the limit contour processes $\zeta^\delta$. We start by showing that the point processes of the rescaled truncated reactant forests converge. We then show that the limit of the point processes for the truncated reactant is equal in distribution to the point processes associated with the contours solving the $(A^\delta, \mathcal{D}(A^\delta))$-martingale problem.

Throughout we assume that $\{\tilde{\eta}^{\text{tot},n}; n \in \mathbb{N}\}$ and $X$ are given on one probability space such that (4.1) holds. Let

$$R^{0,1} := \inf\{t \geq 0 : \xi^{\text{tot}}_t = 0\}.$$

The first lemma describes the distribution of the reactant particle point process (prior to any rescaling). Recall that $\eta^{\text{tot}}$ equals in distribution $\tilde{\eta}^{\text{tot},1}$ by definition.

LEMMA 8.1 (Reactant point processes prior to rescaling). *For all $t \in (0, R^{0,1})$, the reactant point process $\Xi^t$ is a simple point process $\{(i, \tau_i); i = 1, \ldots, \xi^{\text{tot}}_t - 1\}$, where given $\xi^{\text{tot}}_t$, the $\{\tau_i; i = 1, \ldots, \xi^{\text{tot}}_t - 1\}$ are independent and identically distributed with*

$$\mathbb{P}(\tau_1 \in dh) = \frac{\eta^{\text{tot}}_h}{(1 + \int_h^t \eta^{\text{tot}}_s \, ds)^2} \cdot \frac{1 + \int_0^t \eta^{\text{tot}}_s \, ds}{\int_0^t \eta^{\text{tot}}_s \, ds} \cdot \mathbf{1}_{(0,t)}(h) \, dh.$$

PROOF. By Lemma 6.1, the pair $(C, V)$ has the strong Markov property. Hence the excursions below $t$ of $C$ are independent and identically distributed.

Notice that an excursion below $t$ equals in distribution the contour process of a reactant tree run from time $t$ backward, that is, an upside down tree whose branching rates are governed by the process $(\eta^{\text{tot}}_{t-s})_{0 \leq s \leq t}$, and which is conditioned to have height less than $t$. The maximal depth of the excursion is equivalent in distribution to the extinction time of the reactant process run with branching rates $(\eta^{\text{tot}}_{t-s})_{0 \leq s \leq t}$ conditioned to go extinct by time $t$.

For a critical binary branching process with time-inhomogeneous rates $(\lambda_s)_{s \geq 0}$ the population size process $(N_s)_{s \geq 0}$ satisfies (see, e.g., [15], XVII.10.11)

$$\mathbb{P}\{N_t = 0\} = \frac{\int_0^t \lambda_s \, ds}{1 + \int_0^t \lambda_s \, ds}. \tag{8.1}$$

With the special choice $\lambda_s := \eta^{\text{tot}}_{t-s}$, we get by explicit calculation

$$\mathbb{P}\{\tau_1 \geq h\} = \mathbb{P}(N_{t-h} = 0 | N_t = 0) = \frac{\mathbb{P}\{N_{t-h} = 0\}}{\mathbb{P}\{N_t = 0\}}$$

$$= \frac{\int_h^t \eta^{\text{tot}}_s \, ds}{1 + \int_h^t \eta^{\text{tot}}_s \, ds} \frac{1 + \int_0^t \eta^{\text{tot}}_s \, ds}{\int_0^t \eta^{\text{tot}}_s \, ds}.$$

Differentiating this expression with respect to $h$ gives our claim.  □

We can now proceed with the proof of Theorem 3.

PROOF OF THEOREM 3.   For all $n \in \mathbb{N}$, recall that $\tilde{T}^{0,n}$ from (3.1) is the first time the rescaled mass process of the catalyst hits 0, and let

$$\tilde{R}^{0,n} := \inf\{t \geq 0 : \tilde{\xi}_t^{\mathrm{tot},n} = 0\}.$$

We first examine the distribution of the number of trees in a truncated reactant forest and the genealogical distances for a fixed $n$.

STEP 1 (The point process in the $n$th approximation step).   Choose a sequence $(t_n)_{n \to \infty}$ and $t > 0$ such that $t_n \underset{n \to \infty}{\longrightarrow} t$. Recall that $\Pi^t$ and $\Xi^t$ are the genealogical point processes and $\tilde{\Pi}^{t_n,n}$ and $\tilde{\Xi}^{t_n,n}$ are their rescaled versions from Definitions 3.5 and 3.8. In order to determine the number of trees in the forest, notice that given that $t_n \leq T^{0,n} \wedge R^{0,n}$ and given the value of $\tilde{\xi}_{t_n}^{\mathrm{tot},n}$, by independence of the line segments stated in Lemma 6.1, we have that the number $\mathbf{t}_n := \tilde{\Xi}^{t_n,n}(\{\frac{1}{n}, \ldots, \tilde{\xi}^{\mathrm{tot},n} - \frac{1}{n}\} \times \{0\})$ of points at height 0 satisfies

$$(\mathbf{t}_n; \tilde{\eta}^{\mathrm{tot},n}) \overset{d}{=} \mathrm{Bin}(n\tilde{\xi}_{t_n}^{\mathrm{tot},n} - 1, \mathbb{P}\{\tilde{\tau}^n \geq t_n\}),$$

where $\tilde{\tau}^n$ is defined as the extinction time of $\tilde{\xi}^{\mathrm{tot},n}$ with branching rates $(n \times \tilde{\eta}_{t-s}^{\mathrm{tot}})_{0 \leq s \leq t}$ started with one individual, and $\mathrm{Bin}(n, p)$ denotes a random variable which is binomially distributed with parameters $n$ and $p$.

Fix an arbitrary $k_n \in \{1, \ldots, n\tilde{\xi}_{t_n}^{\mathrm{tot},n} - 1\}$. Let $\kappa_n := \tilde{\Xi}^{t_n,n}(\{\frac{1}{n}, \ldots, \frac{k_n}{n}\} \times \{0\})$ be the random number of points between $\frac{1}{n}$ and $\frac{k_n}{n}$ at height 0. Since we choose the $k_n$ first among the $n\tilde{\xi}_t^{\mathrm{tot},n} - 1$ many excursions we get from the i.i.d. structure that

$$(\kappa_n; \tilde{\eta}^{\mathrm{tot},n}) \overset{d}{=} \mathrm{Bin}(k_n, \mathbb{P}\{\tilde{\tau}^n \geq t_n\}).$$

Furthermore, given the values of $(\kappa_n; \tilde{\eta}^{\mathrm{tot},n})$,

$$
\begin{aligned}
&\left(\tilde{\Xi}^{t_n,n}\left(\left\{\frac{1}{n}, \ldots, \frac{k_n}{n}\right\} \times (0, t_n - h_n]\right); \tilde{\eta}^{\mathrm{tot},n}\right) \\
&\qquad \overset{d}{=} \mathrm{Bin}(k_n - (\kappa_n; \tilde{\eta}^{\mathrm{tot},n}), \mathbb{P}(\tilde{\tau}^n \geq h_n | \tilde{\tau}^n < t_n)).
\end{aligned}
$$
(8.2)

STEP 2 (The limit point process).   We now evaluate the asymptotics as $n \to \infty$. Pick a sequence $(k_n)_{n \in \mathbb{N}}$ on $\mathbb{N}$ such that $\frac{k_n}{n} \underset{n \to \infty}{\longrightarrow} uY_t$, for some $u \in [0, 1]$. Since (8.1) implies, given the catalyst total mass processes $\tilde{\eta}^{\mathrm{tot},n}$ and $X$, that for $(h_n) \in [0, t]$ such that $h_n \to h$ for some $h \in [0, t]$,

$$n\mathbb{P}\{\tilde{\tau}^n \geq h_n\} = \frac{n}{1 + n \int_{t_n - h_n}^{t_n} \tilde{\eta}_s^{\mathrm{tot},n} \, ds} \underset{n \to \infty}{\longrightarrow} \left(\int_{t-h}^t X_s \, ds\right)^{-1},$$

it follows, by the usual Poisson approximation of the Binomial distribution, that there exists a point process $\aleph^t$ such that

$$(8.3) \qquad (\kappa_n; \tilde{\eta}^{\mathrm{tot},n}) \underset{n \to \infty}{\Longrightarrow} (\aleph^t([0, uY_t] \times \{0\}); X) \overset{d}{=} \mathrm{Pois}\left( \frac{uY_t}{\int_0^t X_s\, ds} \right).$$

Next fix $h \in (0, t)$, and let $(h_n)$ be a sequence in $(0, t_n)$ such that $h_n \to h$. Note that the first argument of the right-hand side in (8.2) is by (8.3) asymptotically equivalent to $un Y_t$. The second argument is of order $n^{-1}$. Namely by Lemma 8.1,

$$n\mathbb{P}\{\tilde{\tau}^n \geq h_n | \tilde{\tau}^n < t_n\}$$

$$= \frac{n\mathbb{P}\{\tilde{\tau}^n \geq h_n\} - n\mathbb{P}\{\tilde{\tau}^n \geq t_n\}}{\mathbb{P}\{\tilde{\tau}^n < t_n\}}$$

$$= \frac{n^2 \int_0^{t_n - h_n} \tilde{\eta}_s^{\mathrm{tot},n}\, ds}{(1 + n \int_{t_n - h_n}^{t_n} \tilde{\eta}_s^{\mathrm{tot},n}\, ds)(1 + n \int_0^{t_n} \tilde{\eta}_s^{\mathrm{tot},n}\, ds)} \frac{1 + n \int_0^{t_n} \tilde{\eta}_s^{\mathrm{tot},n}\, ds}{n \int_0^{t_n} \tilde{\eta}_s^{\mathrm{tot},n}\, ds}$$

$$\underset{n \to \infty}{\longrightarrow} \frac{1}{\int_{t-h}^t X_s\, ds} - \frac{1}{\int_0^t X_s\, ds}.$$

Hence, once more using the Poisson approximation,

$$\left( \tilde{\Xi}^{t_n, n}\left( \left\{ \frac{1}{n}, \ldots, \frac{k_n}{n} \right\} \times (0, t_n - h_n] \right); \tilde{\eta}^{\mathrm{tot},n} \right)$$

$$(8.4) \qquad \underset{n \to \infty}{\Longrightarrow} (\tilde{\Xi}^t([0, uY_t] \times (0, t - h]); X)$$

$$\overset{d}{=} \mathrm{Pois}\left( uY_t \left( \frac{1}{\int_{t-h}^t X_s\, ds} - \frac{1}{\int_0^t X_s\, ds} \right) \right).$$

Since the limiting intensity measure is absolutely continuous with respect to the Lebesgue measure, (8.4) implies the weak convergence of the point processes, and hence Theorem 3 is proved. $\quad\square$

We next show that the reactant limit point process can be derived from the unique solutions of the $(A^\delta, \mathcal{D}(A^\delta))$-martingale problem. We begin with the proof of part (a), that is, the existence of a local-time process related to $\zeta^\delta$.

PROOF OF PROPOSITION 4.3. (a) For that purpose we rely on the standard characterization (see, e.g., Section VII.3 in [25]) of a one-dimensional diffusion by its scale function $s : [0, \infty) \to \mathbb{R}$ and speed measure $m$, which is based on its construction from a Brownian motion.

Suppose $\zeta$ is an $\mathbb{R}$-valued diffusion with scale function $s$ and speed measure $m$, such that $s$ is differentiable and $m$ is absolutely continuous with respect to the Lebesgue measure. Let $m(dx) = m'(x)\,dx$. Then the generator of the diffusion $\zeta$ is given by

$$(8.5) \qquad Af(x) = \frac{1}{2}\frac{d}{dm}\left(\frac{d}{ds}f(x)\right)$$

with domain given by

$$\mathcal{D}(A) := \left\{ f : f = C + \int_0^{\cdot} s'(x)g(x)\,dx, \text{ where } \frac{g'}{m'} \in C^0([0,\infty)) \right\}.$$

Then the process $B := (B_u)_{u\geq 0}$ defined by $B_u := s(\zeta_u), u \geq 0$, is a diffusion whose scale function is the identity. Recall that since the scale function is nondecreasing it has a well-defined inverse

$$s^{-1}(y) := \inf\{x : s(x) \geq y\}.$$

Then the speed measure $m_B$ and the diffusion coefficient $\sigma_B^2$ of $B$ are

$$m_B(du) = \frac{m'(s^{-1}(u))}{s'(s^{-1}(u))}\,ds, \qquad \sigma_B^2(s^{-1}(u)) = \frac{s'(s^{-1}(u))}{m'(s^{-1}(u))}$$

(compare also with Lemma 6.6).

Let $\beta := (\beta_u)_{u\geq 0}$ be a standard Brownian motion with local-time process $\ell^t(\beta)$ at level $t$ and $\gamma$ be the random-time change defined by

$$(8.6) \qquad \gamma(u) := \int_0^u \frac{m'(s^{-1}(\beta_v))}{s'(s^{-1}(\beta_v))}\,dv.$$

Martingale theory implies that the function $\gamma$ is such that for its inverse $\gamma^{-1}(v) := \{u : \gamma(u) \geq v\}$ we have

$$(B_u)_{u\geq 0} \stackrel{d}{=} (\beta_{\gamma^{-1}(u)})_{u\geq 0}.$$

Consequently, we have a representation of the diffusion $\zeta$ by a distributionally equivalent process

$$(8.7) \qquad (\zeta_u)_{u\geq 0} \stackrel{d}{=} (s^{-1}(\beta_{\gamma^{-1}(u)}))_{u\geq 0}. \qquad\qquad \square$$

We next derive a relation between the local times of $\beta$ and $\zeta$.

LEMMA 8.2 (Excursion process of a diffusion). *Let $\beta$ be the standard Brownian motion, and $\zeta$ be the diffusion with scale function $s$ and speed measure density $m'$. Let $\ell^t(\zeta)$ be the local-time process of $\zeta$ at level $t$. Let $\varepsilon_\ell^-$ denote the downward excursions of $\zeta$ from level $t$ indexed by the local time $\ell^t(\zeta)$, and let $\pi^{\zeta,t} := \{(\ell, \inf(\varepsilon_\ell^-))\}$ be the point process of their depths.*

(i) *The local-time processes* $\ell^t(\zeta) := (\ell^t_u(\zeta))_{u \geq 0}$ *of the diffusion* $\zeta$ *at the level* $t > 0$ *and* $\ell^t(\beta) := (\ell^t_u(\beta))_{u \geq 0}$ *of Brownian motion* $\beta$ *at level* $t > 0$ *satisfy*

$$\ell^t_u(\zeta) \overset{d}{=} \frac{1}{s'(t)} \ell^{s(t)}_{\gamma^{-1}(u)}(\beta).$$

(ii) *The point process* $\pi^{\zeta,t}$ *has the following law: the marks on* $[0,1] \times \{t\}$ *separating contributions from distinct excursions of* $\zeta$ *above* 0 *that reach level* $t$ *are points of a Poisson process with rate*

$$s'(t)(s(t) - s(0))^{-1}.$$

*Each such excursion path of* $\zeta$ *above* 0 *contributes a set of downward excursions* $\varepsilon^-_\ell$ *of* $\zeta$ *below level* $t$ *whose depths are points of a Poisson point process with intensity measure*

$$\aleph^{\zeta,t} = d\ell \otimes n^{\zeta,t}(dh)$$

*and*

(8.8) $$n^{\zeta,t}(dh) = \frac{s'(h)\,dh}{(s(t) - s(h))^2} s'(t).$$

PROOF. (i) By definition,

$$\begin{aligned}
\ell^t_u(\zeta) &= \lim_{\varepsilon \to 0} \frac{1}{2\varepsilon} \int_0^u \mathbf{1}_{\{\zeta_v \in (t-\varepsilon, t+\varepsilon)\}} \, d\langle \zeta, \zeta \rangle_v \\
&= \lim_{\varepsilon \to 0} \frac{1}{2\varepsilon} \int_0^u \mathbf{1}_{\{B_v \in (s(t-\varepsilon), s(t+\varepsilon))\}} \left( \frac{1}{s'(s^{-1}(B_v))} \right)^2 d\langle B, B \rangle_v \\
&= \frac{1}{s'(t)} \ell^{s(t)}_u(B).
\end{aligned}$$

Recall that $\gamma$ is the random-time change from (8.6). We have

$$\gamma(u) = \int_0^u \frac{m'(s^{-1}(\beta_v))}{s'(s^{-1}(\beta_v))} \, dv = \int \ell^t_u(\beta) \frac{m'(s^{-1}(t))}{s'(s^{-1}(t))} \, dt.$$

Hence,

$$\begin{aligned}
\ell^t_u(B) &= \lim_{\varepsilon \to 0} \frac{1}{2\varepsilon} \int_0^u \mathbf{1}_{\{B_v \in (t-\varepsilon, t+\varepsilon)\}} \, d\langle B, B \rangle_v \\
&= \lim_{\varepsilon \to 0} \frac{1}{2\varepsilon} \int_0^{\gamma^{-1}(u)} \mathbf{1}_{\{\beta_v \in (t-\varepsilon, t+\varepsilon)\}} \frac{m'(s^{-1}(\beta_v))}{s'(s^{-1}(\beta_v))} \frac{s'(s^{-1}(\beta_v))}{m'(s^{-1}(\beta_v))} \, d\langle \beta, \beta \rangle_v \\
&= \ell^t_{\gamma^{-1}(u)}(\beta)
\end{aligned}$$

and consequently

(8.9) $$\ell^t_u(\zeta) = \frac{1}{s'(t)} \ell^{s(t)}_{\gamma^{-1}(u)}(\beta).$$

(ii) By Itô's excursion theory $\pi^{\zeta,t}$ is Poisson with an intensity measure (see, e.g., [24], Vol. 2, VII.47)

$$\aleph^{\zeta,t} = d\ell \otimes n^{\zeta,t}(\inf(\varepsilon_\ell^-) \in d(t-h))$$

for some $\sigma$-finite measure $n^{\zeta,t}$ on $[0,t]$. Every excursion $\varepsilon_\ell^-$ of $\zeta$ below $t$ corresponds to a time changed excursion of $\beta$ below $s(t)$. Since the number of such excursions in the two processes is related by (8.9), we have

$$n^{\zeta,t}\{\inf(\varepsilon_\ell^-) > t - h\} = s'(t)n^{\beta,s(t)}\{\inf(\varepsilon_\ell^-) > s(t) - s(h)\}$$
$$= s'(t)(s(t) - s(h))^{-1}.$$

Now, each excursion of $\zeta$ above $0$ that reaches height $t$ has local time at level $t$ equal to $\ell^t_{\tau_{t,0}}(\zeta)$ where $\tau_{t,0} = \inf\{u \geq 0 : \zeta_u = 0,$ if $\zeta_0 = t\}$. Excursion theory also gives that the local time of a Markov process at level $t$ run until the first hitting time of level $0$ is distributed as an exponential variable with parameter $n^{\zeta,t}\{\sup(\varepsilon_\ell^-) > t\} = s'(t)(s(t) - s(0))^{-1}$. Hence the marks on level $t$ separating the contribution to the point process $\pi^{\zeta,t}$ from each such distinct excursion form a Poisson point process on $[0,1] \times \{t\}$ with this rate. $\square$

We will next use Lemma 8.2 to give the proof of Proposition 4.3(b).

PROOF OF PROPOSITION 4.3. (b) Fix $\delta > 0$, and let $\zeta^\delta$ be the solution of the $(A^\delta, \mathcal{D}(A^\delta))$-martingale problem.

Let $\pi^{\zeta^\delta,t}$ be the point process of depths of downward excursions of $\zeta^\delta$ below some fixed level $t$, indexed by the local-time process $\ell^t(\zeta^\delta)$. Recall from (6.11) and (8.5) that $s$ is the scale function and $m$ is the speed measure. By (8.8) with (6.11) $\pi^{\zeta^\delta,t}$ is Poisson with intensity measure

$$\aleph^{\zeta^\delta,t}(du \otimes dh) = d\ell \otimes n^{\zeta,t} = [d\ell^t_u(\zeta^\delta)] \otimes X_t\left[\frac{X_h\, dh}{(\int_h^t X_u\, du)^2}\right].$$

The marks separating contribution from distinct excursions reaching level $t$, in other words from different trees, form a Poisson point process on $[0,1] \times \{t\}$ with rate $X_t(\int_0^t X_s\, ds)^{-1}$.

Let $u^* := \inf\{u \geq 0 : \ell^0_u(2|\beta|) = 1\}$ be the first time the local time of the standard reflected Brownian motion $2|\beta|$ at level $0$ reaches $1$.

Using Lemma 8.2 we have for all $t < \tau^\delta$,

$$(\ell^t_{u^*}(|\zeta^\delta|)X_t)_{0 \leq t < \tau^\delta} \stackrel{d}{=} (\ell^{s(t)}_{\gamma^{-1}(u^*)}(2|\beta|))_{0 \leq t < \tau}.$$

Recall that the second Ray–Knight theorem implies that the local-time process $\ell^t_{u^*}(2|\beta|)$ as a function of the level $t$ until the fixed time $u^*$ is distributed as a Feller diffusion, that is,

$$(\ell^t_{u^*}(2|\beta|))_{t \geq 0} \stackrel{d}{=} (Z_t)_{t \geq 0},$$

where $Z$ satisfies the first part of (3.3), that is,

$$dZ_t = \sqrt{Z_t}\,dW_t^Z \qquad \text{with } Z_0 = 1.$$

If $\gamma^* := \gamma(u^*)$, then (8.9) implies that

$$(\ell_{\gamma^*}^t(|\zeta^\delta|)X_t)_{0<t<\tau^\delta} \stackrel{d}{=} (\ell_{u^*}^{s(t)}(2|\beta|))_{0<t<\tau^\delta} \stackrel{d}{=} (Y_t)_{0<t<\tau^\delta},$$

where $Y$ satisfies the second part of (3.3) with another independent set of Brownian motions, that is,

$$dY_t = \sqrt{s'(t)Y_t}\,dW_t^Y = \sqrt{Z_tY_t}\,dW_t^Y \qquad \text{with } Y_0 = 1.$$

Hence the local-time process of $\zeta^\delta$ at level $t$ equals in distribution $(Y_t)_{t\geq 0}$ from (3.3). The claim of Proposition 4.3 follows. $\quad\square$

**9. Proof for the joint law (Theorem 4).** In this section we show the convergence for the *joint* law of the catalyst and reactant forests.

PROOF OF THEOREM 4. In order to show (4.8) it is enough to show the joint convergence of the pair consisting of the catalyst total mass process and the reactant forest, that is,

$$(9.1) \qquad\qquad (\tilde{\eta}^{\mathrm{tot},n}, \tilde{\xi}^{\mathrm{for},n}) \underset{n\to\infty}{\Longrightarrow} (X, Y^{\mathrm{for}}).$$

This follows since given the values of the catalyst mass process the reactant forest and the catalyst forest are independent.

So far we have shown that $\tilde{\eta}^{\mathrm{tot},n} \underset{n\to\infty}{\Longrightarrow} X$ and that under assumption (4.1), $(\tilde{\xi}^{\mathrm{for},n};\tilde{\eta}^{\mathrm{tot},n}) \underset{n\to\infty}{\Longrightarrow} (Y^{\mathrm{for}};X)$. Hence, in order to show (9.1), it suffices to show that the law of the reactant forest depends continuously on the catalyst mass path. For that purpose recall from Theorem 2 that for all $\delta > 0$, $(Q_\delta(Y^{\mathrm{for}});X) \stackrel{d}{=} (\mathcal{T}(\zeta^\delta);X)$. Since the map $\mathcal{T}:C_{\mathbb{R}^+}^0[0,\infty): \to \mathbb{T}^{\mathrm{root}}$ is continuous, the claim follows by observing from (8.7) that $\zeta^\delta$ is continuous in $X := (X_t)_{t\in(0,\tau^0)}$. $\quad\square$

**10. Proof of differences between reactant and "classical" forest.** In this section we prove the results on the comparison of the reactant limit forest with the classical forest.

PROOF OF PROPOSITION 4.4. We have specified the distribution of both $Z^{\mathrm{for}}$ and $Y^{\mathrm{for}}$ using collections of point processes $\pi^{\zeta,t}$ and $\pi^{\beta,t}$. For each tree in the forest $Y^{\mathrm{for}}$ the distribution of $\pi^{\zeta,t}$ is given by the intensity measure $\aleph^{\zeta,t}$ from (4.7), and for a standard Brownian tree in $Z^{\mathrm{for}}$ distribution of $\pi^{\beta,t}$ is given by the intensity measure $\aleph^{\beta,t}$ in (3.8). It is easy to check that for each $t \in (0, \tau^0]$,

(10.1) $$((\mathbf{s_t}(\aleph^{\zeta,t}); X))_{t \geq 0} \overset{d}{=} ((\aleph^{\beta,t}))_{t \geq 0},$$

where $\mathbf{s_t}$ applies to the second coordinate of the bivariate measure $\aleph^{\zeta,t}$. Since the genealogical point process records the distances between points at the same level in the forest, (10.1) implies the claim. $\square$

PROOF OF PROPOSITION 4.5.

$$\mathbb{E}\left[\int_{(\partial Q_t(Y^{\text{for}}))^2} (\mu^{t,Y})^{\otimes 2}(d(u, u')) \mathbf{1}\{d_{Y^{\text{for}}}(u, u') = 2t\}\right]$$

$$= \int_0^1 du_1 \int_0^1 du_2 \, \mathbb{P}\{\pi^t([u_1 Y_t, u_2 Y_t] \times \{0\}) \neq \varnothing\}$$

$$= 1 - \int_0^1 du_1 \int_0^1 du_2 \, \mathbb{E}\left[\exp\left\{-\frac{Y_t|u_1 - u_2|}{\int_0^t X_s \, ds}\right\}\right]$$

$$= 1 - \int_0^1 du_1 \int_0^1 du_2 \, \mathbb{E}\left[\exp\left\{-\left(\int_0^t X_s \, ds\right)^{-1} \frac{|u_1 - u_2|}{1 + |u_1 - u_2|}\right\}\right],$$

where we have used that a branching diffusion $\tilde{Y}_t$ with time-dependent branching rates $b(t)$ satisfies

$$\mathbb{E}[e^{-\lambda \tilde{Y}_t} | \tilde{Y}_0 = y] = \exp\left\{-\frac{y \cdot \lambda}{1 + \lambda \int_0^t b(u) \, du}\right\}.$$

By assumption $\mathbb{E}[(\int_0^t ds \, X_s)^{-1}] = zt^{-1}$. Therefore, Jensen inequality yields

$$\mathbb{E}\left[\exp-\left\{\left(\int_0^t X_s \, ds\right)^{-1} \frac{|u_1 - u_2|}{1 + |u_1 - u_2|}\right\}\right]$$

$$\geq \exp-\left\{\mathbb{E}\left[\left(\int_0^t X_s \, ds\right)^{-1}\right] \frac{|u_1 - u_2|}{1 + |u_1 - u_2|}\right\}$$

$$= \exp-\left\{\frac{z}{t} \cdot \frac{|u_1 - u_2|}{1 + |u_1 - u_2|}\right\}.$$

Hence

$$\mathbb{E}\left[\int_{(\partial Q_t(Y^{\text{for}}))^2} (\mu^{t,Y})^{\otimes 2}(d(u, u')) \mathbf{1}\{d_{Y^{\text{for}}}(u, u') = 2t\}\right]$$

$$\leq 1 - \int_0^1 du_1 \int_0^1 du_2 \, \exp\left\{-\frac{z}{t} \cdot \frac{|u_1 - u_2|}{1 + |u_1 - u_2|}\right\}$$

$$= \mathbb{E}\left[\int_{(\partial Q_t(Z^{\text{for}}))^2} (\mu^{t,Z})^{\otimes 2}(d(u, u')) \mathbf{1}\{d_{Z^{\text{for}}}(u, u') = 2t\}\right]. \qquad \square$$

We prepare the proof of Corollary 4.7 with the following lemma. Assume that we are given an excursion $e$ with $L(e) = \inf\{t > 0 : e(t) = 0\} < \infty$. Recall from (2.1) the map $\mathcal{T}$ which sends $e$ to a rooted ordered forest $\mathcal{T}(e) = (T_e, d_{T_e}, [0], \leq_{T_e})$ associated with $e$.

LEMMA 10.1. *The $\mathbf{l}^2$-length of a forest $\mathcal{T}(e)$ associated with an excursion $e \in C^0_{\mathbb{R}^+}[0, \infty)$ with $L(e) = \inf\{t > 0 : e(t) = 0\} < \infty$ equals the quadratic variation of $e$*

$$(10.2) \qquad \chi_e := P - \lim_{\varepsilon \to 0} \sup_{\substack{0 \leq t_0 < t_1 < \cdots < t_n \leq L(e), \\ \max\{|t_{i+1} - t_i| \leq \varepsilon\}}} \left\{ \sum_{i=1}^{n} |e(t_{i+1}) - e(t_i)|^2 \right\}.$$

PROOF. Assume that we are given an excursion $e$ with $L(e) = \inf\{t > 0 : e(t) = 0\} < \infty$. Assume that $T_e$ has a well-defined $\mathbf{l}^2$-length. Since edges of positive length do not contribute to the $\mathbf{l}^2$-lengths of $T_e$, (4.9) holds even if we further restrict the finite subsets $T'_n$ of $T_e$ to those with $\inf\{\varepsilon > 0 : (T'_n)^\varepsilon = T_e^{\text{branch}}\} \to 0$, as $n \to \infty$. Here $T_e^{\text{branch}}$ denotes the set of all inner branch points of $T_e$, where $x \in T_e$ is an inner branch point if $T_e \setminus \{x\}$ decomposes $T_e$ into more than 2 connected components. Likewise, monotone regions in the excursion $e$ do not contribute to the quadratic variation of $e$, so in (10.2) we may restrict the set of time points to only the finite subsets of the set of times $0 \leq t_0^{\min} < \cdots < t_i^{\min} < \cdots \leq L(e)$ such that $e(t_i^{\min})$ is a local minimum. Since $e$ is continuous it can have only a countable number of local minima on an interval $[0, L(e)]$. Now we observe that the local minima of $e$ precisely correspond to the branch points of $T_e$, so

$$\chi_e = P - \lim_{\varepsilon \to 0} \sup_{\max\{|t_{i+1}^{\min} - t_i^{\min}| \leq \varepsilon\}} \left\{ \sum |e(t_{i+1}^{\min}) - e(t_i^{\min})|^2 \right\}$$

$$= P - \lim_{\varepsilon \to 0} \sup_{\max\{|t_{i+1}^{\min} - t_i^{\min}| \leq \varepsilon\}} \left\{ \sum \left| e(t_{i+1}^{\min}) + e(t_i^{\min}) - 2 \min_{s \in [t_i^{\min}, t_{i+1}^{\min}]} e(s) \right|^2 \right\}$$

$$= P - \lim_{n \to \infty} \sum_{a,b \in T'_n} (d(a,b))^2$$

$$= \mathbf{l}^2(\mathcal{T}(e)). \qquad \square$$

PROOF OF COROLLARY 4.7. Lemma 10.1 states in particular that if $\zeta$ is assumed to be an excursion associated with the reactant forest then $\chi_\zeta = \infty$ by Lemma 4.6. Hence $\zeta$ cannot be a diffusion path and the claim follows. $\square$

PROOF OF LEMMA 4.6. By its definition (4.9), $\mathbf{l}^2(Z^{\text{for}})$ equals in distribution the quadratic variation of a reflected Brownian motion $2|\beta|_{\cdot \wedge \ell^{-1}(1)}$ which is known to be $\ell(2|\beta|)^{-1}(1) < \infty$, almost surely.

As for the reactant limit tree, on the event $\{\rho^0 < \tau^0\}$, $\mathbf{l}^2(Y^{\text{for}})$ equals in distribution for suitably small positive $\delta$ the quadratic variation of the well-posed solution $\zeta^\delta$ of the $(A^\delta, \mathcal{D}(A^\delta))$-martingale problem stopped at the first time the local time of $\zeta^\delta$ at 0 reaches 1 $(\ell^0(\zeta^\delta))^{-1}(1)$. This quadratic variation is finite almost surely and remains finite as $\delta \to 0$.

On the other hand, on $\{\rho^0 \geq \tau^0\}$ we show that the quadratic variation of the contours diverges as $\delta \to 0$. In fact, by (10.2),

$$\mathbf{l}^2(Q_{\tau^\delta}(Y^{\text{for}}; X)) = \langle \zeta^\delta, \zeta^\delta \rangle_{(\ell^0(\zeta^\delta))^{-1}(1)}.$$

Recall from (8.7) that, given $X$, the law of $\zeta^\delta$ can be expressed in terms of a function of a time changed Brownian motion $\beta$:

$$(\zeta_u^\delta)_{u \geq 0} \overset{d}{=} (s^{-1}(B_u))_{u \geq 0}, \qquad (B_u)_{u \geq 0} \overset{d}{=} (\beta_{\gamma^{-1}(u)})_{u \geq 0}.$$

Since $\gamma(u) = \int_0^u \frac{1}{s'(s^{-1}(\beta_v))} dv$ we have $d\langle B, B \rangle_u = s'(s^{-1}(B_u)) du$, and

$$d\langle \zeta^\delta, \zeta^\delta \rangle_u \overset{d}{=} \left( \frac{1}{s'(s^{-1}(B_u))} \right)^2 d\langle B, B \rangle_u = \frac{1}{s'(s^{-1}(B_u))} du = \frac{1}{X_{\zeta_u^\delta}} du.$$

Hence, by the occupation time formula

$$\mathbf{l}^2(Q_{\tau^\delta}(Y^{\text{for}}; X)) = \int_0^{\ell^0(\zeta^\delta)^{-1}(1)} \frac{1}{X_{\zeta_v^\delta}} dv \overset{d}{=} \int_0^{s(\tau^\delta)} \frac{\ell(\beta)_{\gamma(\ell^0(\beta)^{-1}(1))}^t}{X_{s^{-1}(t)}^2} dt$$

$$= \int_0^{\tau^\delta} \frac{\ell(\beta)_{\gamma(\ell^0(\beta)^{-1}(1))}^{s(v)}}{X_v} dv$$

by changing variables. We next note that $\inf\{v \geq 0 : \ell(\beta)_{\gamma(\ell^0(\beta)^{-1}(1))}^{s(v)} = 0\} \overset{d}{=} \rho^0 > \tau^0$, so we will have $c = \inf_{v \in [0, \tau^0]} \ell(\beta)_{\gamma(\ell^0(\beta)^{-1}(1))}^{s(v)} dv > 0$ and

$$\mathbf{l}^2(Q_{\tau^\delta}(Y^{\text{for}}; X)) \geq c \int_0^{\tau^0} \frac{1}{X_v} dv \underset{\delta \to 0}{\to} \infty$$

using the known fact that $\int_0^{\tau^0} \frac{1}{X_v} dv = \infty$ (see, e.g., Lemma 3.1 in [5]). $\square$

**Acknowledgment.** The authors thank the referee for extremely careful reading of the paper. Her or his suggestions and comments have led to many improvements in the paper.

## REFERENCES

[1] ALDOUS, D. (1991). The continuum random tree. I. *Ann. Probab.* **19** 1–28. MR1085326

[2] ALDOUS, D. (1991). The continuum random tree. II. An overview. In *Stochastic Analysis (Durham, 1990). London Mathematical Society Lecture Note Series* **167** 23–70. Cambridge Univ. Press, Cambridge. MR1166406

[3] ALDOUS, D. (1993). The continuum random tree. III. *Ann. Probab.* **21** 248–289. MR1207226

[4] ALDOUS, D. J. (1999). Deterministic and stochastic models for coalescence (aggregation and coagulation): A review of the mean-field theory for probabilists. *Bernoulli* **5** 3–48. MR1673235

[5] ATHREYA, S. and WINTER, A. (2005). Spatial coupling of neutral measure-valued population models. *Stochastic Process. Appl.* **115** 891–906. MR2134483

[6] CHISWELL, I. (2001). *Introduction to Λ-Trees.* World Scientific, River Edge, NJ. MR1851337

[7] DONNELLY, P. and KURTZ, T. G. (1999). Particle representations for measure-valued population models. *Ann. Probab.* **27** 166–205. MR1681126

[8] DRESS, A. W. M. (1984). Trees, tight extensions of metric spaces, and the cohomological dimension of certain groups: A note on combinatorial properties of metric spaces. *Adv. in Math.* **53** 321–402. MR753872

[9] DRESS, A. W. M., MOULTON, V. and TERHALLE, W. F. (1996). *T*-theorie. *Europ. J. Combinatorics* **17** 161–175.

[10] DRESS, A. W. M. and TERHALLE, W. F. (1996). The real tree. *Adv. in Math.* **120** 283–301. MR1397084

[11] DUQUESNE, T. and LE GALL, J.-F. (2002). Random trees, Lévy processes and spatial branching processes. *Astérisque* **281** vi–147. MR1954248

[12] ETHIER, S. N. and KURTZ, T. G. (1986). *Markov Processes: Characterization and Convergence.* Wiley, New York. MR838085

[13] EVANS, S. N., PITMAN, J. and WINTER, A. (2006). Rayleigh processes, real trees, and root growth with re-grafting. *Probab. Theory Related Fields* **134** 81–126. MR2221786

[14] EVANS, S. N. and WINTER, A. (2006). Subtree prune and re-graft: A reversible real-tree valued Markov chain. *Ann. Probab.* To appear. **34** 918–961.

[15] FELLER, W. (1968). *An Introduction to Probability Theory and Its Applications. Vol. I,* 3rd ed. Wiley, New York. MR0228020

[16] GREVEN, A., KLENKE, A. and WAKOLBINGER, A. (1999). The longtime behavior of branching random walk in a catalytic medium. *Electron. J. Probab.* **4** 80. MR1690316

[17] KURTZ, T. G. (1992). Averaging for martingale problems and stochastic approximation. In *Applied Stochastic Analysis (New Brunswick, NJ, 1991). Lecture Notes in Control and Information Sciences* **177** 186–209. Springer, Berlin. MR1169928

[18] LAMPERTI, J. (1967). The limit of a sequence of branching processes. *Z. Wahrsch. Verw. Gebiete* **7** 271–288. MR0217893

[19] LE GALL, J.-F. (1999). *Spatial Branching Processes, Random Snakes and Partial Differential Equations.* Birkhäuser, Basel. MR1714707

[20] PENSSEL, C. (2003). Interacting Feller diffusions in catalytic media. Ph.D. thesis, Institute of Math., Erlangen, Germany.

[21] PITMAN, J. (2006). *Combinatorial Stochastic Processes. Lecture Notes in Mathematics* **1875**. Springer, Berlin. MR2245368

[22] POPOVIC, L. (2004). Asymptotic genealogy of a critical branching process. *Ann. Appl. Probab.* **14** 2120–2148. MR2100386

[23] PROTTER, P. E. (1977). On the existence, uniqueness, convergence and explosions of solutions of systems of stochastic integral equations. *Ann. Probab.* **5** 243–261. MR0431380

[24] ROGERS, L. C. G. and WILLIAMS, D. (1987). *Diffusions, Markov Processes, and Martingales: Itô Calculus. Wiley Series in Probability and Mathematical Statistics: Probability and Mathematical Statistics.* **2**. Wiley, New York. MR921238

[25] REVUZ, D. and YOR, M. (1999). *Continuous Martingales and Brownian Motion*, 3rd ed. *Grundlehren der Mathematischen Wissenschaften [Fundamental Principles of Mathematical Sciences]* **293**. Springer, Berlin. MR1725357

[26] STROOCK, D. W. and VARADHAN, S. R. S. (1969). Diffusion processes with continuous coefficients. I. *Comm. Pure Appl. Math.* **22** 345–400. MR0253426

[27] TERHALLE, W. F. (1997). *R*-trees and symmetric differences of sets. *European J. Combin.* **18** 825–833. MR1478827

A. GREVEN
A. WINTER
MATHEMATISCHES INSTITUT
UNIVERSITY ERLANGEN–NURNBERG
BISMARCKSTR. 1 1/2
91054 ERLANGEN
GERMANY
E-MAIL: greven@mi.uni-erlangen.de
          winter@mi.uni-erlangen.de

L. POPOVIC
DEPARTMENT OF MATHEMATICS AND STATISTICS
CONCORDIA UNIVERSITY
MONTREAL QC H3G 1M8
CANADA
E-MAIL: lpopovic@mathstat.concordia.ca